\documentclass[11pt]{amsart}
\usepackage{latexsym,amssymb}
\usepackage{graphics}
\usepackage{url}
\usepackage{harvard}
\usepackage{eulervm}
\usepackage{epsfig,subfigure}
\usepackage{tikz}
\newcommand*\circled[1]{\tikz[baseline=(char.base)]{
            \node[shape=circle,draw,inner sep=1pt] (char) {#1};}}
\usepackage{graphicx}

\newtheoremstyle{normal}
{2ex}               
{3ex}               
{}                  
{}                  
{\bfseries} 
{}                  
{2pt}   
{\thmname{#1}\thmnumber{ #2.} \thmnote{(#3)}}

\theoremstyle{normal}
\newtheorem{definition}{Definition}
\newtheorem{remark}[definition]{Remark}
\newtheorem{example}{Example}[section]

\newtheorem{lemma}{Lemma}
\newtheorem{assumption}{Assumption}
\newtheorem{theorem}{Theorem}
\newtheorem{corollary}{Corollary}
\newtheorem{proposition}{Proposition}

\renewenvironment{remark}{\noindent {\bf Remark.\ }} {\hfill{ \rule{0mm}{0mm}}}

\DeclareMathOperator{\Var}{Var}

\DeclareMathOperator{\Cov}{Cov}

\def\FF{\mathcal F}
\def\E{\mathbb E}
\def \eps {\epsilon}
\title{Neyman's $C(\alpha)$ test for unobserved heterogeneity }
\author{Jiaying Gu \\ University of Illinois at Urbana-Champaign}
\date{October 4, 2014}
\thanks{Department of Economics, University of Illinois at Urbana-Champaign, 214 David Kinley Hall, 1407 W. Gregory Dr., Urbana, Illinois 61801, MC-707, USA. Tel: +1-267-994-1519. Fax: +1-217-244-6571. Email Address: gu17@illinois.edu. I would like to thank Roger Koenker for his continued support and encouragement. I would also like to thank Andreas Hagemann, Marc Hallin, Keisuke Hirano, Stanislav Volgushev, two anonymous referees and the participants at the Midwest Econometrics Group meeting 2012 and the Boneyard Conference 2013 at the University of Illinois for valuable comments and useful discussion. I gratefully acknowledge financial support from NSF grant
SES-11-53548 and the Paul Boltz summer Fellowship. All errors are my own.  }                                           

\begin{document}
\bibliographystyle{econometrica}

\begin{abstract}
A unified framework is proposed for tests of unobserved heterogeneity in parametric statistic models based on Neyman's $C(\alpha)$ approach. Such tests are irregular in the sense that the first order derivative of the log likelihood with respect to the heterogeneity parameter is identically zero, and consequently the conventional Fisher information about the parameter is zero. Nevertheless, local asymptotic optimality of the $C(\alpha)$ tests can be established via LeCam's differentiability in quadratic mean and the limit experiment approach. This leads to local alternatives of order $n^{-1/4}$. The scalar case result is already familiar from existing literature and we extend it to the multi-dimensional case. The new framework reveals that certain regularity conditions commonly employed in earlier developments are unnecessary, i.e. the symmetry or third moment condition imposed on the heterogeneity distribution. Additionally, the limit experiment for the multi-dimensional case suggests modifications on existing tests for slope heterogeneity in cross sectional and panel data models that lead to power improvement. Since the $C(\alpha)$ framework is not restricted to the parametric model and the test statistics do not depend on the particular choice of the heterogeneity distribution, it is useful for a broad range of applications for testing parametric heterogeneity. 
\end{abstract}

\maketitle

\section{Introduction}

\possessivecite{Neyman59} $C(\alpha)$ test can be viewed as a generalization of \possessivecite{Rao48} score test in the presence of nuisance parameters and thus provides a unified framework for parametric statistical inference. We will see that many of the existing tests for neglected parameter heterogeneity can also be formulated as $C(\alpha)$ tests and share common features. However, for these tests the usual score function is identically zero under the null hypothesis, and conventional Fisher information is thus zero. Fortunately, in these cases the second derivative of the log likelihood is non-degenerate and approximations based on it can be used to form a modified version of LeCam's differentiability in quadratic mean (DQM) condition. Local asymptotic normality (LAN) theory, then leads to local asymptotic optimality results for the $C(\alpha)$ test in such settings under local alternatives of order $n^{-1/4}$.\\

We find that LeCam's limit experiment perspective is very useful in analyzing tests for neglected heterogeneity especially in the multi-dimensional setting. It allows us to first develop optimal test statistics for the Gaussian limit and then extend them to the corresponding asymptotic $C(\alpha)$ test. The one-sided nature of the limit experiment reveals that we require the mixture of $\chi^2$ asymptotics which leads to power improvement compared to the conventional $\chi^2$ type test. This finding is relevant to the Information Matrix test and some of the recent applications to slope heterogeneity testing in panel data models. \\

We focus initially on the case of a scalar heterogeneity parameter. Although some of the results are already familiar in the literature, the use of the LeCam framework is new and it leads to a set of less restrictive assumptions and sheds light on why reparameterization leads to unnecessary conditions employed in previous literature. Discussing the scalar case in the LeCam framework also facilitates the extension to multivariate settings which is described at the end of Section 2 and is the major contribution of the paper. In Section 3 we consider four different examples. In the first example, the $C(\alpha)$ tests for parameter heterogeneity in Poisson regression model under two slightly different alternative specifications lead to tests introduced in \citeasnoun{Lee86}. The second example considers testing for slope heterogeneity in cross sectional linear regression models; the $C(\alpha)$ test in this setting shares many features of the \citeasnoun{Breusch79} LM test, but the positivity constraints revealed via the limit experiment suggest a modification that leads to a power gain. We then illustrate an example using the $C(\alpha)$ test to jointly test for heterogenous location and scale parameters in Gaussian panel data models. Lastly, we compare the $C(\alpha)$ test for slope heterogeneity in panel data model to the test considered in \citeasnoun{PY08}. For a wide range of $N$ and $T$, the $C(\alpha)$ test, since it pays explicit attention to the positivity constraints under the alternative, enjoys a power improvement.  \\

The $C(\alpha)$ test for heterogeneity formulated in this paper is very similar to the setup used in some previous development. In a seminal paper, \citeasnoun{Chesher84} points out the score test for unobserved parametric heterogeneity is identical to \possessivecite{White82} Information Matrix (IM) test. \citeasnoun{Cox83} obtains similar results using a more general mixture model. These papers can be viewed as important extensions of a somewhat neglected example on testing for parameter heterogeneity in Poisson models in \citeasnoun{NeymanScott}. \citeasnoun{Moran73} investigates the asymptotic behavior of these score tests.  However, as we will show in Section 4, the parameterization adopted in \citeasnoun{Moran73} and also in \citeasnoun{Chesher84} requires some unnecessary additional assumptions, i.e. the zero third moment or symmetry of the heterogeneity distribution, even though it delivers the same score function formed based on the second derivative as the $C(\alpha)$ test constructed here. The explanation is that the likelihood under their parameterization obtains the same expansion of the likelihood under the $C(\alpha)$ test parameterization only if symmetry holds. In addition, even though the score function is the same, the positivity constraints lead to a different decision rule with mixture of $\chi^2$ asymptotics in contrast to the conventional $\chi^2$ test for the IM test. Furthermore, there are situations where the $C(\alpha)$ test for unobserved heterogeneity is no longer identical to the IM test and we illustrate some conditions for equivalence to hold in Section 4. Lastly, a Monte Carlo simulation is carried out to evaluate the power performance for various examples.

\section{The $C(\alpha)$ test for unobserved parameter heterogeneity}

\citeasnoun{Neyman59} introduces the $C(\alpha)$ test with the consideration that hypotheses testing problems in applied research often involve several nuisance parameters. In these composite testing problems, most powerful tests do not exist, motivating search for an optimal test procedure that yields the highest power among the class of tests obtaining the same size. Neyman's locally asymptotically optimality result for the $C(\alpha)$ test employs regularity conditions inherited from the conditions used by \citeasnoun{Cramer} for showing consistency of MLE and some further restrictions on the testing function to allow for replacing the unknown nuisance parameters by its $\sqrt{n}$-consistent estimators. It is the confluence of these Cram\'er conditions and the maintained significance level $\alpha$ that gives the name to the $C(\alpha)$ test. 

\subsection{$C(\alpha)$ test in regular cases} \label{2.1}
In regular cases, where all the score functions with respect to parameters in the model are non-degenerate and the Fisher information matrix is non-singular, the $C(\alpha)$ test is constructed as follows. Suppose we have $X_{1}, \dots, X_{n}$ as i.i.d. random variables with density $p(x; \xi, \theta)$ where $\theta$ are nuisance parameters belonging to $\Theta \subset \mathbb{R}^{p}$ and $\xi$ are parameters under test that belong to $\Xi \subset \mathbb{R}^{q}$. For densities satisfing the regularity conditions (\citeasnoun[Definition 3]{Neyman59}), we consider testing the hypothesis $H_{0}: \xi = \xi_{0}$ against $H_{a}: \xi \in \Xi \setminus \{ \xi_{0} \} $ while nuisance parameters $\theta \in \Theta$ are left unspecified. We define the conventional score functions as
\[
C_{\xi, n} = \frac{1}{\sqrt{n}} \sum_{i=1}^{n} \nabla_{\xi} \log p(X_{i} ; \xi, \theta) |_{\xi=\xi_{0}}
\]
\[
C_{\theta, n} = \frac{1}{\sqrt{n}} \sum_{i=1}^{n} \nabla_{\theta} \log p(X_{i} ; \xi, \theta) |_{\xi=\xi_{0}} 
\]
and denote the corresponding Fisher information matrix as,
\[
I = \left ( 
\begin{array}{cc}
	I_{\xi \xi} & I_{\xi \theta} \\
	I_{\theta \xi} & I_{\theta \theta} 
\end{array} 
\right ).
\]

Since nuisance parameters $\theta$ are left unspecified by $H_{0}$, \citeasnoun{Neyman59} shows that for the test statistic to have the same asymptotic behavior when we replace the nuisance parameters $\theta$ by any $\sqrt{n}$-consistent estimator $\hat \theta_{n}$, it is necessary and sufficient for the test statistics to be orthogonal to $C_{\theta, n}$. For example, the "residual" score, which constitutes the vector of projecting $C_{\xi, n}$ onto the space spanned by the score vector $C_{\theta, n}$, denoted by 
\[
g_{n}(\theta) = C_{\xi, n} - I_{\xi \theta} I_{\theta \theta}^{-1} C_{\theta, n}, 
\]
provides such a test function with variance $I_{\xi . \theta} \equiv I_{\xi \xi}-I_{\xi \theta} I_{\theta \theta}^{-1} I_{\theta \xi}$. Given a $\sqrt{n}$-consistent estimator $\hat \theta_{n}$ for $\theta$, the $C(\alpha)$ test
\[
T_{n} (\hat \theta_{n}) = g_{n}(\hat \theta_{n}) ^\top I_{\xi . \theta}^{-1} g_{n} (\hat \theta_{n})
\]
is then asymptotically $\chi^{2}_{q}$ under $H_{0}$ and is optimal for local alternatives of the form $\xi_{n} = \xi_{0} + \delta / \sqrt{n}$. When $\hat \theta_{n}$ is the restricted maximum likelihood estimator of $\theta$, $C_{\theta, n}$ is zero and the $C(\alpha)$ test reduces to Rao's score test. The component $I_{\xi \theta} I_{\theta \theta}^{-1} I_{\theta \xi}$ subtracted from the information $I_{\xi \xi}$ for $\xi$ measures the amount of information lost due to not knowing the nuisance parameters (see e.g. \citeasnoun{Bickeletal}, section 2.4). \\

\subsection{Testing for unobserved parameter heterogeneity} \label{2.2}
The $C(\alpha)$ test for unobserved heterogeneity is usually formulated under a random parameter model. Following \citeasnoun{NeymanScott} we will focus initially on testing homogeneity of a scalar parameter against the alternative that the parameter is random. Consider having i.i.d. random variables $X_{1}, \dots, X_{n}$, with each $X_{i}$ having density function $p(x; \lambda_{i})$. Heterogeneity of the model is introduced by regarding the individual specific $\lambda_{i}$ as a random parameter of the form, 
\[
\lambda_{i} = \lambda_{0} + \tau \xi U_{i},
\]
where the unobserved $U_{i}$'s are independent random variables with common distribution function, $F$, satisfying moment conditions $\mathbb{E}(U)=0$, $\mathbb{V}(U)=1$. The parameter $\tau$ is a known finite scale parameter, which allows us to rescale the variance for $U$ to be unity. It is not restrictive to assume $\tau$ known, as we will see later that $\tau$ does not enter the test statistics. It is cancelled out when the test function is studentized by its standard deviation. The hypothesis we would like to test is $H_{0}: \xi =0$, which implies $\lambda_{i} = \lambda_{0} $ for all $i$'s. The alternative hypothesis is $H_{a}: \xi \neq 0$. \\

Under the above setup, the standard $C(\alpha)$ test described in Section \ref{2.1} breaks down because the score function for $\xi$ for each individual observation $x_{i}$, defined as the first order logarithmic derivative of the density function with respect to $\xi$, is identically zero under the null, hence the Fisher information is also zero, 
\[
\frac{\partial }{\partial \xi} \log \int p(x_{i}; \lambda_{0} + \tau \xi u) dF(u) \left |_{\xi=0} \right. = \tau \int udF(u) \frac{p^{\prime} (x_{i}; \lambda_{0})}{p(x_{i} ; \lambda_{0})} =0.
\] 
However, in circumstances like this, we can compute the second-order derivative, denoted as $s_i(\lambda_0)$ below, 
\[
s_{i}(\lambda_{0}):=\frac{\partial^{2}}{\partial \xi^{2}} \log \int p(x_{i} ; \lambda_{0} + \tau \xi u) dF(u)\left |_{\xi=0} \right.  = \tau^{2} \int u^{2}dF(u) \frac{p^{\prime \prime}(x_{i}; \lambda_{0}) }{p(x_{i}; \lambda_{0})} = \tau^{2} \frac{p^{\prime \prime}(x_{i}; \lambda_{0})}{p(x_{i}; \lambda_{0})}.
\]
\\
The normed sum of these independent second-order derivatives, $s(\lambda_{0})=\frac{1}{\sqrt{n}} \sum_{i} s_{i}(\lambda_{0})$, can be shown to be asymptotically normally distributed with mean zero and variance $\mathbb{E}(s_{1}^{2}(\lambda_{0}))$ under $H_{0}$ by the central limit theorem and by noticing that $\E(p^{\prime \prime}(x_{i}; \lambda_{0})/p(x_{i}; \lambda_{0}))=0$ as a consequence of differentiating $\int p(x;\lambda)dx=1$ as a function of $\lambda$ twice. This leads to a close analogy with the classical theorem, in which $s(\lambda_{0})$ acts as the score function and the variance $\mathbb{E}(s_{1}^{2}(\lambda_{0}))$ plays the role of the Fisher information in the irregular setting considered here.\\

In regular cases, score tests exploit the fact that if the null hypothesis is false, the gradient of the log likelihood should not be close to zero. Clearly this fails in the irregular case, because no matter how data is generated, the gradient is always zero. It is natural then to make use of the curvature information provided by the second-order derivative for inference. If the null is false, one expects the second-order derivative to be positive. We will see that this second-order score function plays the essential role of constructing the $C(\alpha)$ test for unobserved heterogeneity. The positivity condition also anticipates that the $C(\alpha)$ test will be one-sided. The goal of the remaining part of this section is to show that the optimality of the $C(\alpha)$ test, as in the regular case, is still preserved under this irregularity and its asymptotic theory, although different from the regular cases in certain perspectives, still takes a simple form.

\subsection{Asymptotic optimality of the $C(\alpha)$ test for parameter heterogeneity}
Under the irregularity discussed above, in order to establish the optimality of the test statistics based on the second-order score function, one could consider modifying the Cram\'er type regularity conditions in \citeasnoun[Definition 3]{Neyman59}, requiring the density function to be five times differentiable pointwise and impose a Lipschitz condition on the fifth order derivative with respect to the parameter under test. The main motivation is to obtain a quadratic approximation of the log likelihood ratio using the second-order score function through a higher order Taylor expansion. To be more specific, using the example in Section \ref{2.2} as an illustration, for local alternatives $\lambda_{i} = \lambda_{0} + \tau \xi_{n} U_{i}$, with $\xi_{n}$ be a sequence that converges to zero at certain rate, we have the following Taylor expansion of the log likelihood ratio,
\[
\begin{array}{cc}
\Lambda_{n} = \sum_{i} \log \frac{p(x_{i}; \lambda_{i})}{p(x_{I}; \lambda_{0})} = \frac{\xi_{n}^{2} \tau^{2}}{2!} \mathbb{E}(U^{2})\sum_{i} s_{i}(\lambda_{0}) + \frac{\xi_{n}^{3}\tau^{3}}{3!} \mathbb{E}(U^{3}) \sum_{i} \frac{\nabla_{\lambda}^{3}p(x_{i}; \lambda_{0})}{p(x_{i}; \lambda_{0})} \\
+ \frac{\xi_{n}^{4} \tau^{4}}{4!} \left [ \mathbb{E}(U^{4}) \sum_{i} \frac{\nabla_{\lambda}^{4}p(x_{i}; \lambda_{0})}{p(x_{i}; \lambda_{0})} - 3 \mathbb{E}(U^{2})^{2} \sum_{i} s_{i}^{2}(\lambda_{0})\right ] + o_{P}(1).
\end{array}
\]
Let $\xi_{n}$ be of order $n^{-1/4}$ and provided the third and fourth moments of $U$ are finite in addition to the zero mean and unit variance assumption, we obtain a quadratic approximation of the log-likelihood. More details of such regularity conditions can be found in \citeasnoun{Rotnitzky}, in which they consider the maximum likelihood estimation of $\xi$ in the irregular cases in a very general context. \citeasnoun[Chapter 4]{Lindsay95} also has a brief discussion of this. \\

An alternative formulation, rooted in LeCam's local asymptotic normality (LAN) theory, can be based on 
his differentiability in quadratic mean (DQM) condition. The latter condition is less stringent in regular cases: while Cram\'er conditions assume the density to be three times differentiable and impose a Lipschitz condition on the third order derivative, the DQM condition only requires first order differentiability and the derivative to be square integrable in $\mathcal{L}_{2}$ space. \citeasnoun{Pollard97} provides a nice discussion of the DQM condition in these regular cases. This is the new approach we take for analyzing the asymptotic behavior of the $C(\alpha)$ test for heterogeneity. We will show below that by modifying the DQM condition slightly, we can obtain the local asymptotic normality of the log-likelihood ratio and establish the asymptotic optimality of the $C(\alpha)$ test for the irregular cases under assumptions much weaker than those suggested by the classical Neyman's approach. One prominent example for which the classical conditions fail while the DQM conditions are satisfied is the double exponential location model with $p_{\theta}(x) = f(x-\theta)$ and $f(x) = \frac{1}{2} \exp (-|x|)$. For this model, the density function $f$ is not differentiable at 0 but it satisfies the DQM condition. We would thus have no difficulty constructing a test for homogeneity in the location parameter for this model under the LeCam type conditions.\\

Suppose we have a random sample $(X_{1}, \dots, X_{n})$ with density function $p(x; \xi, \theta)$ with respect to some measure $\mu$. The joint distribution of this i.i.d. random sample will be denoted as $P_{n, \xi, \theta}$, which is the product of $n$ copies of the marginal distribution $P(x; \xi, \theta)$. 

\begin{assumption} 
The density function $p$ satisfies the following conditions:
\begin{enumerate}
\item $\xi_{0}$ is an interior point of $\Xi$
\item For all $\theta \in \Theta \subset \mathbb{R}^{p}$ and $\xi \in \Xi \subset \mathbb{R}$, the density is twice continuously differentiable 
	with respect to $\xi$ and once continuously differentiable with respect to $\theta$ for $\mu$-almost all $x$. 
\item Denoting the first two derivatives of the density with respect to $\xi$ 
	evaluated under the null as $\nabla_\xi p(x; \xi_{0}, \theta)$ and 
	$\nabla_\xi^2 p(x; \xi_{0}, \theta)$, we have 
	$\mathbb{P} \left (\nabla_\xi p(x;  \xi_{0}, \theta)=0 \right)=1$ and 
	$\mathbb{P} \left(\nabla_\xi^2 p(x; \xi_{0}, \theta) \neq 0 \right)>0$ for all $\theta \in \Theta \subset \mathbb{R}^{p}$.
\item Denoting the derivative of the density with respect to $\theta$ evaluated 
	under the null as $\nabla_{\theta}p(x; \xi_{0}, \theta)$, for any 
	$p$-dimensional vector $a$, $\mathbb{P} \left(\nabla_\xi^2 p(x; \xi_{0}, \theta) 
	\neq a^\top \nabla_{\theta}p(x; \xi_{0}, \theta) \right)>0$. 
\end{enumerate}
\end{assumption}

\begin{remark}
Here $\xi$ is the parameter under test and $\theta$ is the vector of nuisance parameters. The list of regularity conditions in Assumption 1 tailors the standard conditions for a regular $C(\alpha)$ test to the heterogeneity test we consider here. In particular, condition (3) reflects the irregularity of these tests that the first order logarithmic derivative with respect to $\xi$ vanishes but the second-order derivative is non-vanishing. Condition (2) secures existence of the respective derivatives. Condition (4) rules out the case where there is a perfect linear relationship between the second-order score for $\xi$ and the score for $\theta$. It ensures the new Fisher information thus defined to be non-singular and the $C(\alpha)$ test statistics to be non-degenerate.  \\
\end{remark}

Under Assumption 1, we can now define the modified DQM condition that is crucial for establishing the local asymptotic normality of the model. \\

\begin{definition}
The density $p(x; \xi, \theta)$ satisfies the modified differentiability in 
quadratic mean condition at $(\xi_{0}, \theta)$ if there exists a 
vector $v(x)=(v_{\xi}(x), v_{\theta}^\top(x))^\top \in \mathcal{L}_{2}(\mu)$ such that as 
$(\xi_{n}, \theta_{n}) \to (\xi_{0}, \theta)$, 
\begin{displaymath}
\int |\sqrt{p(x; \xi_{n}, \theta_{n})}-\sqrt{p(x; \xi_{0}, \theta)} 
- h_{n}^\top v(x)|^{2} d\mu(x)
= o( ||h_{n}||^{2})
\end{displaymath}
where $h_{n}=((\xi_{n}-\xi_{0})^{2}, (\theta_{n}-\theta)^\top)^\top$. Here $|| \cdot ||$  denotes the Euclidean norm and $\mathcal{L}_{2}(\mu)$ denotes the $\mathcal{L}_{2}$ space of square integrable functions with respect to measure $\mu$. \\
\indent Furthermore, let $\beta(h_{n})$ be the mass of the part of $p(x; \xi_{n}, \theta_{n})$ that is $p(x; \xi_{0}, \theta)$-singular, then as $(\xi_{n}, \theta_{n}) \to (\xi_{0}, \theta)$, 
\[
\frac{\beta(h_{n})}{|| h_{n} ||^{2}} \to 0
\]
\end{definition}

Usually the vector $v(x)$ contains derivatives of the square root of density $\sqrt{ p(x; \xi_{n}, \theta_{n})}$ with respect to each parameter evaluated under their null value. Definition 1 modifies the classical DQM condition such that whenever the first order derivative is identically zero for certain parameters, it is differentiated again until it is nonvanishing. The corresponding terms in $h_{n}$ also need to be raised to the same power. For the heterogeneity test, the score function with respect to $\xi$ is of second order and its associated term in $h_{n}$ is hence quadratic. This further implies that the contiguous alternatives must be $O(n^{-1/4})$. For the following theorems, we will thus focus on the sequence of local models on $(X_{1}, \dots, X_{n})$ with joint distribution $P_{n, \xi_{n}, \theta_{n}}$ in which $\xi_{n} =\xi_{0} + \delta_{1} n^{-1/4}$ and $\theta_{n} = \theta + \delta_{2} n^{-1/2}$.  \\

\begin{theorem} \label{theorem1}

Suppose $(X_{1}, \dots, X_{n})$ are i.i.d. random variables with joint distribution $P_{n, \xi_{n}, \theta_{n}}$ and the density satisfies Assumption 1 and the modified DQM condition with 

\begin{displaymath}
v(x)=(v_{\xi}(x), v_{\theta}^\top(x))^\top=\left (\frac{1}{4}  
\frac{\nabla_\xi^2 p(x ; \xi_{0}, \theta)}{\sqrt{p(x; \xi_{0}, \theta)}}  
\mathbb{I} _{[p(x; \xi_{0}, \theta)>0]}, 
\frac{1}{2} \frac{\nabla_{\theta} p(x ; \xi_{0},\theta)^\top}{\sqrt{p(x ; \xi_{0}, \theta)} }  
\mathbb{I} _{[p(x; \xi_{0},\theta)>0]} \right )^\top, 
\end{displaymath}
then for fixed $\delta_{1}$ and $\delta_{2}$, the log-likelihood ratio has the following quadratic approximation under the null:
\begin{displaymath}
\Lambda_{n}=\log \frac{dP_{n, \xi_{n}, \theta_{n}}}{dP_{n,\xi_{0}, \theta}} = t^\top S_{n} -\frac{1}{2} t ^\top J t + o_{P}(1)
\end{displaymath}
where $t=(\delta_{1}^{2}, \delta_{2}^\top)^\top$, 
\[
S_{n}=(S_{\xi, n}, S_{\theta, n}^\top)^\top=\left ( \frac{2}{\sqrt{n}} \sum_{i} \frac{v_{\xi}(x_{i})}{\sqrt{p(x_{i} ; \xi_{0}, \theta)}}, \frac{2}{\sqrt{n}} \sum_{i} \frac{v_{\theta}^\top(x_{i})}{\sqrt{p(x_{i} ; \xi_{0}, \theta)}} \right )^\top
\]
and \\
\[
J=4 \int (v v^\top) d\mu(x)=\begin{pmatrix}
\mathbb{E}(S_{\xi,n}^{2}) & \Cov(S_{\xi,n},S_{\theta,n}^\top) \\ 
\Cov (S_{\xi,n}, S_{\theta,n})& \mathbb{E}(S_{\theta,n}S_{\theta,n}^\top)
\end{pmatrix}
\equiv \begin{pmatrix}
J_{\xi \xi} & J_{\xi \theta}\\ 
J_{\theta \xi} & J_{\theta \theta}
\end{pmatrix}.
\]\\
\end{theorem}

\begin{corollary} \label{cor1}
With $S_{n}$ and $J$ defined as in Theorem 1, we have 
\[
S_{n} \overset{P_{n, \xi_{0}, \theta}}{\leadsto} \mathcal{N}(0, J),
\]
and hence the sequence of models $P_{n, \xi_{n}, \theta_{n}}$ is locally asymptotically normal (LAN) at $(\xi_{0}, \theta)$ with $S_{n}$ being interpreted as the score vector and $J$ as the associated Fisher information matrix. Furthermore, $P_{n, \xi_{n}, \theta_{n}}$ is mutually contiguous to $P_{n, \xi_{0}, \theta}$. 

\end{corollary}

Theorem 1 shows that under Assumption 1, the modified DQM condition is sufficient for obtaining a quadratic approximation of the log-likelihood ratio for the sequence of local models in the $n^{-1/4}$ neighborhood of the null value $\xi_{0}$ and the $n^{-1/2}$ neighborhood of the nuisance parameter $\theta$. The joint normality of the vector $S_{n}$, as established in Corollary 1, further indicates the LAN property of this sequence of models. It is important to note that the vector $S_{n}$, in which the degenerately zero first-order score function for $\xi$ is replaced by the corresponding second-order derivative of the log-likelihood, acts as the score vector in this irregular case. Naturally, $J$ has the interpretation of the Fisher information matrix. Under Assumption 1, since we rule out perfect dependence between $S_{\xi,n}$ and $S_{\theta,n}$ in condition (4), $J$ is non-singular. \\

Having established the LAN property of this sequence of local models, we can now make use of \possessivecite{LeCam72} limit experiment theory to show that the $C(\alpha)$ test is locally asymptotically optimal in the scalar case. \\

Following the definitions given in \citeasnoun{LeCam72} and \citeasnoun{vdv98}, an experiment $\mathcal{E}$ indexed by a parameter set $H$ is a collection of probability measures $\{ P_{h}: h \in H \}$ on the sample space $(\mathcal{X}, \mathcal{A})$. A sequence of experiments $\mathcal{E}_{n} = (\mathcal{X}_{n}, \mathcal{A}_{n}, P_{n, h}: h \in H)$ is said to converge to a limit experiment $\mathcal{E}=(\mathcal{X}, \mathcal{A}, P_{h}: h \in H)$ if the likelihood ratio process for $\mathcal{E}_{n}$, $\frac{dP_{n, h}}{dP_{n, h_{0}}}(X_{n})$, converges in distribution to the likelihood ratio of the limit experiment, $\frac{dP_{h}}{dP_{h_{0}}}(X)$, for $h$ in every finite subset $I \subset H$ and for every null value $h_{0} \in H$ . A common feature is that many sequences of experiments produce a Gaussian limit experiment. One important example is that for i.i.d. sample from a smooth parametric model with distribution $P_{\vartheta}$, if the sequence of the local model $P_{n, \vartheta_{n}}$ in which $\vartheta_{n} = \vartheta_{0} + r_{n} \delta $ with $r_{n}$ as the appropriate norming rate is locally asymptotically normal, then it has a Gaussian shift experiment as its limit. \\

The advantage of establishing the limit experiment is several fold. First, the limit experiment is often easier to analyze than the original sequence of models. Second, the limit experiment provides a bound for the optimal estimation (in terms of lower bound on the asymptotic variance) or testing procedure (in terms of upper bound on the asymptotic power) one could achieve in the original model. Third, by the asymptotic representation theory (\citeasnoun[Chapter 9]{vdv98}), any sequence of statistics that converges in the original experiment can be matched in the limit experiment and they share identical asymptotic behavior. We will show in particular that the $C(\alpha)$ test statistics is matched with the optimal testing procedure in the Gaussian shift limit experiment, hence establishing its optimality. We first focus on the scalar case, leaving the multi-dimensional case to a separate discussion. \\

\begin{theorem}\label{theorem2}
Let $\mathcal{E}_{n}$ be a sequence of experiments based on i.i.d. random variables $(X_{1}, \dots, X_{n})$ with joint distribution $P_{n, \xi_{n}, \theta_{n}}$ on the sample space $(\mathcal{X}_{n}, \mathcal{A}_{n})$. We further index the sequence of experiment by $t =(\delta_{1}^{2}, \delta_{2}^\top)^\top \in \mathbb{R}_{+} \times \mathbb{R}^{p}$. The log-likelihood ratio of the sequence of models satisfies, 
\[
\log \left (\frac{dP_{n,\xi_{n}, \theta_{n}}}{dP_{n, \xi_{0}, \theta}} \right )
= t^\top S_{n} -\tfrac{1}{2} t^\top J t + o_{P}(1),
\]
with the score vector $S_{n}$ defined as in Theorem 1 converging in distribution under the null to $\mathcal{N}(0, J)$. Then the sequence of experiments $\mathcal{E}_{n}$ converges to the limit experiment based on observing one sample from $Y= t + v$, where $v \sim \mathcal{N}(0, J^{-1})$. The locally asymptotically optimal statistic for testing $H_{0}: \delta_{1}=0$ vs. $H_{a}: \delta_{1} \neq 0$ is
\[
Z_{n}=(J_{\xi \xi}-J_{\xi \theta}J_{\theta \theta}^{-1}J_{\theta \xi})^{-1/2}(S_{\xi, n}-J_{\xi \theta}J_{\theta \theta}^{-1}S_{\theta,n}).
\] 
\end{theorem}

\begin{corollary} \label{cor2}
Under $H_{0}$, $Z_{n}$ has distribution $\mathcal{N}(0,1)$. Under $H_{a}$, by applying LeCam's third lemma (see e.g. \citeasnoun[Example 6.7]{vdv98}), it follows a shifted normal distribution $\mathcal{N}(\delta_{1}^{2}(J_{\xi \xi}-J_{\xi \theta}J_{\theta \theta}^{-1} J_{\theta \xi})^{1/2},1)$.
\end{corollary}

The optimal test statistic $Z_{n}$ takes the form of a $C(\alpha)$ test. It projects the second-order score $S_{\xi, n}$ for $\xi$ onto the space spanned by the first-order score vector $S_{\theta, n}$ for $\theta$. It is the sequence of statistics from the original experiment that can be matched with the optimal test statistic in the limit Gaussian experiment for inference on $\delta_{1}$, which is the first element in the one sample $Y$.  \\

One common feature of $C(\alpha)$ heterogeneity tests is that the limit distribution under local alternative is always a right-shifted normal distribution even if we have a two-sided alternative hypothesis for $\delta_{1}$. This is not surprising given that the shift parameter corresponding to $\xi$ in the Gaussian limit experiment is a quadratic term $\delta_{1}^{2} \in \mathbb{R}_{+}$. In other words, the best inference procedure one could possibly achieve in the limit experiment is for $\delta_{1}^{2}$. We lose the sign information on $\delta_{1}$, and the asymptotically optimal test, if rejects the null, fails to distinguish whether the deviation is from the left or from the right (this phenomenon is also emphasized in \citeasnoun{Rotnitzky}). Let $Y = (Y_{1}, Y_{2})^\top$ where the partition is such that  $Y_{1}$ is a scalar and $Y_{2} \in \mathbb{R}^{p}$ as in Theorem \ref{theorem2}. In the Gaussian limit experiment based on the one sample from $Y \sim \mathcal{N}(t, J^{-1})$, the one-sided test, rejecting $H_{0}$ if $Y_{1} \geq \Phi^{-1}(1-\alpha) (J_{\xi \xi} - J_{\xi \theta} J_{\theta \theta}^{-1} J_{\theta \xi})^{-1/2}$, is the uniformly most powerful test. Since the sequence that converges to the rescaled first element $(J_{\xi \xi} - J_{\xi \theta} J_{\theta \theta}^{-1} J_{\theta \xi})^{1/2} Y_{1}$ is exactly $Z_{n}$, it implies that the asymptotic $C(\alpha)$ test rejects $H_{0}$ if $Z_{n} \geq \Phi^{-1}(1-\alpha)$ for any level $\alpha$. Observe that for $\alpha < 0.5$, this is equivalent to rejecting $H_{0}$ if $(0 \vee Z_{n})^{2} > c$, where $c$ is the $(1-\alpha)$-quantile of $\frac{1}{2}\chi^{2}_{0} + \frac{1}{2} \chi^{2}_{1}$ and $\chi^{2}_{0}$ is a degenerate distribution with mass 1 at 0. No solution exists for $c$ if $\alpha >0.5$ although this is of little relevance in practice. We mention the mixture of $\chi^{2}$ asymptotics just to be more cohesive with the multi-dimensional extension later. The weight $1/2$ associated with $\chi^{2}_{0}$ is due to the fact that $Z_{n}$ takes negative values with probability $1/2$ under $H_{0}$. \\

There is another intuitive interpretation of the one-sidedness of the test, as we have already anticipated in Section \ref{2.2}. The $C(\alpha)$ test $Z_{n}$, constructed from the second-order score for $\xi$, exploits information of the curvature of the log-likelihood function. Since at $\xi=\xi_{0}$, the gradient of the log-likelihood function with respect to $\xi$ is always zero, it depends on the sign of the second-order derivative to determine whether the null point is a local maximum or a local minimum. Only positive values of $Z_{n}$ indicates the null point as a local minimum of the log-likelihood function, leading to a rejection of the null hypothesis. As $n \to \infty$, due to normality of $Z_{n}$, only half the time we get the "correct" curvature allowing us to reject the null. In the simulation exercise in Section \ref{4}, we show that paying attention to this one-sided feature in the decision rule gives more power on testing for parameter heterogeneity. \\

For the random parameter model, one could of course also consider a likelihood ratio test as an alternative testing strategy for heterogeneity. Among many others, \citeasnoun{Chenetal01} considers a modified likelihood ratio test for homogeneity in finite mixture models, which is very close to the setup we consider in this paper. They also obtain a mixture of $\chi^{2}$ asymptotics for their likelihood ratio test statistics. Their modified LRT can be viewed as an asymptotically equivalent testing procedure in mixture models to the $C(\alpha)$ test considered here. The latter, however, inheriting the nice feature of the score test, is much easier to compute. Furthermore, the $C(\alpha)$ test statistics does not depend on the specification of $F$ as long as the moment conditions are satisfied. This can be viewed as a merit of the test because it has power for a large class of alternative models. On the other hand, it can also be viewed as its disadvantage because rejecting the hypothesis does not provide information on what plausible alternatives might be. Comparison between the general LR test for mixture models and the $C(\alpha)$ test is considered in \citeasnoun{GKS}. \\

The result established thus far is not confined to the heterogeneity test problem. It is applicable whenever the first-order score for the parameter under test vanishes but the second-order score is non-degenerate. There is another possible scenario for the score test to break down, in which none of the first-order score function is vanishing, but there is linear dependence among them, and thus the Fisher information matrix becomes singular. This is the case discussed in considerable detail in \citeasnoun{LeeChesher}. Models with selection bias and the stochastic production frontier models fall into this class. They propose an extremum test which is based on the determinant of the matrix of the second-order derivatives of the log likelihood function and show the asymptotic optimality of the test. The extremum test can essentially be reformulated, using a reparameterization slightly different from what the authors suggested in the paper (i.e. choose $k$ to be 1 in \citeasnoun[p. 132]{LeeChesher}), to fit into the conditions described in Assumption 1. The similar irregularity also arises in test for symmetry in normal-skew distribution and is investigated in \citeasnoun{Hallin12}. The reparameterization is a Gram-Schmidt orthogonalization in the same spirit of \citeasnoun[Section 4.4]{Rotnitzky}. The $C(\alpha)$ test can then be constructed and asymptotic optimality of the test follows. 

\subsection{Replacing the nuisance parameter by a $\sqrt{n}$-consistent estimator}

Notice that the optimal test statistic $Z_{n}$ we obtained in Theorem 2 is a function of $\theta$, to make the test statistic feasible under unknown nuisance parameters, we need to replace $\theta$ by some estimator $\hat \theta$. In order to ensure that the asymptotics  for the test statistic $Z_{n}$ in Corollary \ref{cor2} is still valid, it suffices to show that $Z_{n}(\hat \theta)- Z_{n}(\theta)=o_{P}(1)$ both under the null and local alternatives. There are various ways to obtain this result. The classical approach taken in \citeasnoun{Neyman59} was to make additional differentiability and bound conditions on the test function $g(x_i, \theta)$, which is defined as 
\[
g(x_i;  \theta)=(J_{\xi \xi}-J_{\xi \theta}J_{\theta \theta}^{-1}J_{\theta \xi})^{-1/2} \Big( \frac{2v_{\xi}(x_i)}{\sqrt{p(x_i; \xi_0, \theta)}}-J_{\xi \theta}J_{\theta \theta}^{-1}\frac{2v_{\theta}(x_i)}{\sqrt{p(x_i; \xi_0, \theta)}}\Big),
\]
such that $Z_{n}(\theta)=\frac{1}{\sqrt{n}} \sum_{i} g(x_{i}, \theta)$. Details of these assumptions can be found in \citeasnoun[Definition 3 (ii) (iii)]{Neyman59} and we will not replicate them here. When the conditions are satisfied, Taylor expansion of $Z_{n}( \hat\theta)$ around $Z_{n}(\theta)$ yields the desired results for $\hat \theta$ being any $\sqrt{n}$-consistent estimator for $\theta$. Neyman's assumptions are rather strong, for example, he requires the density to be three times differentiable with respect to $\theta$ and also moments of the gradient of $g$ with respect to $\theta$ to be continuous. LeCam proposes a discretization trick which works as long as the model satisfies a uniform LAN condition and the $\sqrt{n}$-consistent estimator satisfies an asymptotic discreteness property. The trick is quite standard in one-step estimation problems. Our approach, using more modern probability theory,  is to view the difference $Z_n(\hat \theta) - Z_n(\theta)$ as an empirical process. More precisely, we make the following assumption on the test function $g(x, \theta)$ to establish the equicontinuity of the empirical process. In fact our Assumption \ref{assumption2} below on $g(x, \theta)$ implies the conditions of the Type IV function in \citeasnoun{Andrew94} with $p=2$. We include these results here for completeness. 

\begin{assumption} \label{assumption2}
There exists some $\delta>0$ such that for any $\eta,\eta'\in U_\delta(\theta)$ we have for some $\gamma>0$
\[
|g(x,\eta) - g(x,\eta')| \leq \|\eta-\eta'\|^\gamma H(x)
\] 
for $P_{n,\xi_n,\theta}$-almost all $x$ (for every $n \in \mathbb{N}$) where $H$ is square integrable with respect to $P_{n,\xi_n,\theta}$ for all $n\in \mathbb{N}$, $\sup_n \mathbb{E}_{P_{n,\xi_n,\theta}}H^2(X)<\infty$ and additionally for some $c_n=o(1)$, $n^{1/2}\mathbb{E}_{P_{n,\xi_n,\theta}}[H(X)\mathbb{I}_{\{H(X)>n^{1/2}c_n\}}] = o(1)$. 
\end{assumption}

\begin{theorem} \label{theorem3}
Under Assumption \ref{assumption2}, if $\hat \theta$ is a $\sqrt{n}$-consistent estimator for $\theta$, then
\[
|Z_{n}(\hat \theta)-Z_{n}(\theta)| =o_{P}(1)
\]
\end{theorem}

\subsection{$C(\alpha)$ test for parameter heterogeneity in higher dimensions} \label{2.3.3}

It is of interest to generalize the $C(\alpha)$ tests of unobserved parameter heterogeneity to higher dimensions in the irregular case. For example, in a linear regression model, we may want to jointly test for slope heterogeneity for more than one covariates. When panel data is available, we may want to test for heterogeneity in the slope coefficients in the presence of individual variances, see for example \citeasnoun{PY08}. The main challenge comes from the one-sidedness of the test. Fortunately, the limit experiment turns out to be multivariate Gaussian with location shifts in each coordinate (or in a subset of coordinates) towards the right tail. This naturally requires us to look for optimal tests for deviations of the location parameters of the multivariate Gaussian from zero restrictions to the positive orthant.\\

To be more specific, suppose the limit multivariate Gaussian  experiment has mean vector $(\mu_{1}, \dots, \mu_{q})$, we would like to test $H_{0}: \mu_{i} = 0$ for $i = 1, \dots, q$ against the alternative $H_{a}: \mu_{i} \geq 0$ for $i = 1, \dots, q$ with at least one inequality holds strictly. Unlike in the univariate case where the one-sided test is optimal in the sense of being uniformly most powerful and hence the asymptotic analogue $C(\alpha)$ test obtains the same optimality locally asymptotically,  there exists no uniformly optimal test for the multivariate case. There are two dominant options in the literature. The likelihood ratio test has been studied by many authors. \citeasnoun{Chernoff54} extends the classical Wilks's result on likelihood ratio test (LRT) to cases in which the null value of the parameters under test lie on the boundary of the parameter space. \citeasnoun{Bartholomew}, \citeasnoun{Nuesch} and \citeasnoun{Perlman} among many others consider variants of Gaussian LRT under restricted alternatives. \citeasnoun{Hillier86} provides details for the LRT with dimension equal to three. \citeasnoun{SelfLiang87} give some further examples for LRT with nuisance parameters. Alternatively, \citeasnoun{AbelsonTukey} propose tests based on the idea of maximin contrast and this is further extended by \citeasnoun{SchaafsmanSmid} who introduce the optimality concept of ``most stringent somewhere most powerful'' (MSSMP) test. \\

Neither LRT nor MSSMP test uniformly dominates each other, but both are shown to be substantially more powerful than the usual $\chi^2$ or $F$ test for the multivariate Gaussian case. We construct the $C(\alpha)$ test by extending the LRT via the limit experiment into its local asymptotic version. It allows a direct power comparison to the usual $\chi^2$ test, i.e. the Information Matrix test, which ignores the positivity constraints. It is also closer to the historical development of generalization of the regular $C(\alpha)$ test to the multidimensional case by \citeasnoun{Bulher}, which can be viewed as the asymptotic analogue of the usual $\chi^2$ test in the Gaussian limit experiment for testing $\mu_i = 0$ against non-constrained alternative $\mu_i \neq 0$. Additionally, as we will show, the $C(\alpha)$ test can also be easily adapted if only a subset of the shift parameters are subject to positivity constraints. \\

The LRT statistics for these one-sided test problems in multi-dimensions all obtain a mixture of $\chi^{2}$ with different degrees of freedom as their asymptotic distribution. One disadvantage of the LRT is that the weights of these $\chi^{2}$'s get complicated very quickly as dimension increases in most cases. In contrast, the MSSMP test has a standard normal asymptotics. However, it is hard to adapt the MSSMP test to situations where only a subset of the shift parameters are subject to constraints. \\

We will present in details the joint test for heterogeneity in dimension two as an illustration and comment on the more general case. Suppose again we have i.i.d. random sample $(X_{1}, \dots, X_{n})$ with density $p(x ; \xi, \theta)$. The parameters under test are now $\xi = (\xi_{1}, \xi_{2}) \in \Xi \subset \mathbb{R}^{2}$. They take null value $\xi_{0} = (\xi_{10}, \xi_{20})$ and $\theta \in \Theta \subset \mathbb{R}^{p}$ are the nuisance parameters. For heterogeneity tests in particular, we consider testing for heterogeneity of a vector of parameters, $\lambda_i$, of the model. Under the alternative, they take the form, $\lambda_{ki}=\theta_{k} + \tau \xi_{k} U_{ki}$, for $k=1,2$. Let the covariance matrix for $U_{i} = (U_{1i},U_{2i})$ be $\Omega$. Without loss of generality, we let the diagonal element of $\Omega$ be unity. Under $H_{0}, \xi_{k}=0$, so that $\lambda_{k}$'s are homogenous across individuals taking value $\theta_{k}$. \\

The density function satisfies Assumption 1 such that the first-order score vector for $\xi_{1}$ and $\xi_{2}$ are vanishing due to the zero mean assumption for $U_{i}$ but each elements in the second-order score matrix are non-vanishing. It also satisfies the modified DQM condition so that the model is locally asymptotically normal. Typically the score function for $(\xi_{1}, \xi_{2})$ then consists of all distinct elements in the second-order score matrix. Depending on the assumption on $\Omega$, some of the elements become zero. For example, if $\Omega$ is a diagonal matrix, which implies that $U_{1i}$ is mutually independent to $U_{2i}$, then the off-diagonal terms of the second-order score matrix for $\xi$ are zero. If $\Omega$ has non-zero off-diagonal elements, then the corresponding cross terms in the score matrix are also non-vanishing and need to be included. \\

It is crucial to distinguish the above-mentioned two scenarios, since the diagonal terms in the score matrix correspond to the shift terms in the Gaussian limit experiment that are subject to positivity constraints, while the off-diagnonal terms correspond to shift parameters that can take value over the whole real line. This implies that if $\Omega$ is not diagonal, then the Gaussian limit experiment has only a subset of the shift parameters that have positivity constraints under the alternative. Theorem \ref{theorem4} gives the general theory on constructing the $C(\alpha)$ statistics with the subsequent Corollary \ref{corollary4} discussing the special case if $\Omega$ is diagonal. \\

To proceed, we denote the second-order score vector (all distinct elements in the second order score matrix stacked into a vector) for $(\xi_{1}, \xi_{2})$ as $(S_{\xi_{1}^2,n}, S_{\xi_{2}^2,n}, S_{\xi_{1}\xi_{2}, n})$. The first two corresponds to the diagonal terms and the last the off-diagonal term. Let the first-order score for $\theta$ be $S_{\theta,n }$. More specifically, under regularity conditions, they are $S_{\xi_{k}^2,n}=\frac{1}{2\sqrt{n}} \sum_{i} \frac{\nabla_{\xi_{k} \xi_{k}}^{2}p(x_{i} ;  \xi_{0}, \theta)}{p(x_{i} ; \xi_{0}, \theta)}\mathbb{I}_{[p(x_{i} ; \xi_{0}, \theta)>0]}$ for $k = 1,2$; and $S_{\xi_{1} \xi_{2}, n} = \frac{1}{2 \sqrt{n}} \sum_i \frac{\nabla_{\xi_1 \xi_2}^{2} p (x_i ; \xi_0, \theta)}{p(x_i; \xi_0, \theta)}\mathbb{I}_{[p(x_{i}; \xi_{0}, \theta)>0]}$ and $S_{\theta,n}=\frac{1}{\sqrt{n}}\sum_{i} \frac{\nabla_{\theta} p(x_{i}; \xi_{0}, \theta)^\top}{p(x_{i}; \xi_{0}, \theta)}\mathbb{I}_{[p(x_{i}; \xi_{0}, \theta)>0]}$. Let the associated information matrix be denoted as, $J=\begin{pmatrix}
J_{\xi \xi} &J_{\xi \theta} \\ 
J_{\theta \xi} & J_{\theta \theta}
\end{pmatrix}$, with $J_{\xi \xi}$ being a $3 \times 3$ block matrix. The residual score for $\xi$, similar to the scalar case, is found to be
\[
\tilde S_{\xi,n}=\begin{pmatrix}
\tilde S_{\xi_{1}^2,n} \\
\tilde S_{\xi_{2}^2,n}\\
\tilde S_{\xi_{1}\xi_{2}, n}
\end{pmatrix}:=\begin{pmatrix}
S_{\xi_{1}^2,n}\\ 
S_{\xi_{2}^2,n}\\
S_{\xi_{1}\xi_{2},n}
\end{pmatrix}
 -J_{\xi \theta }J_{\theta \theta}^{-1}S_{\theta,n}
\]
and the covariance matrix for $\tilde S_{\xi, n}$ is $\Sigma = J_{\xi \xi}-J_{\xi \theta}J_{\theta \theta}^{-1}J_{\theta \xi} := \begin{pmatrix} \Sigma_{(11)} & \Sigma_{(12)} \\
\Sigma_{(21)} & \Sigma_{(22)}
\end{pmatrix}$. The partition of $\Sigma$ is such that $\Sigma_{(11)}$ collects covariance terms for the first two elements in $\tilde S_{\xi, n}$. 

\begin{theorem} \label{theorem4}
Let $\upsilon_n$ be the sequence of experiments based on i.i.d. random variable $(X_{1}, \dots, X_{n})$ with joint distribution $P_{n, \xi_{n}, \theta_{n}}$ with $\xi_{n} = (\xi_{10}, \xi_{20}) + (\delta_{1}, \delta_{2})n^{-1/4}$ and $\theta_{n} = \theta +\delta_{3} n^{-1/2} $ on the sample space $(\mathcal{X}_{n}, \mathcal{A}_{n})$. The log-likelihood ratio of the sequence of experiment satisfies, 
\[
\log \left ( \frac{dP_{n, \xi_{n}, \theta_{n}}}{dP_{n, \xi_{0}, \theta}} \right ) = t^\top S_{n} - \frac{1}{2} t ^\top J t + o_{p}(1),
\]
with $S_{n} = (S_{\xi_{1}^2, n}, S_{\xi_{2}^2, n}, S_{\xi_{1}\xi_{2}, n}, S_{\theta, n}^\top)^\top \sim \mathcal{N}(0, J)$. Then the limit experiment of $\upsilon_n$ is based on observing one sample from $Y = t + \nu$ with $t = (\delta_{1}^{2}, \delta_{2}^{2}, 2\delta_{1}\delta_{2}, \delta_{3}^\top)^\top \in \mathbb{R}_{+}^{2} \times \mathbb{R} \times \mathbb{R}^{p}$ and $\nu \sim \mathcal{N}(0, J^{-1})$. We would like to jointly test $H_{0}: \delta_{1}=\delta_{2} =0$ against the alternative $H_{a}: \delta_{1} \neq 0 $ or $\delta_{2} \neq 0$. Let $u_{n} := (u_{1n}, u_{2n})^\top = (\tilde S_{\xi_{1}^2, n}, \tilde S_{\xi_{2}^2,n})^\top - \Sigma_{(12)}\Sigma_{(22)}^{-1} \tilde S_{\xi_{1}\xi_{2},n}$ and let $\Lambda$ be the Cholesky decompositon of $\Sigma_{11.2}:= \Sigma_{(11)} - \Sigma_{(12)}\Sigma_{(22)}^{-1}\Sigma_{(21)}$, that is
\[
\Lambda = \begin{pmatrix}
\sqrt{v_{1}} & 0 \\
\rho \sqrt{v_{2}} & \sqrt {v_{2}} \sqrt{1-\rho^{2}}
\end{pmatrix}
\]
where $\rho$ is the correlation coefficient between $u_{1n}$ and $u_{2n}$ and $v_{1}$ and $v_{2}$ are their respective variances. Define $w_{n}=(w_{1n}, w_{2n})^\top$ as 
\[
w_{n}
\equiv \Lambda^{-1}u_{n} =  \begin{pmatrix} 
u_{1n}/\sqrt{v_{1}}\\
(1-\rho^{2})^{-1/2} ( u_{2n}/\sqrt{v_{2}} - \rho u_{1n} / \sqrt{v_{1}}) 
\end{pmatrix}
\]
and let $w_{3n} := \Sigma_{(22)}^{-1/2} S_{\xi_{1}\xi_{2},n}$. The $C(\alpha)$ test statistic is one of the following four cases:
\[
T_{n} = \begin{cases}
w_{1n}^{2}+w_{2n}^{2} +w_{3n}^{2}& \mbox{if  } w_{1n} \geq \frac{\rho}{\sqrt{1-\rho^{2}}} w_{2n}, w_{2n} \geq 0\\
w_{1n}^{2}  +w_{3n}^{2}& \mbox{if  } w_{2n} \leq 0, w_{1n} \geq 0\\
(\rho w_{1n} + \sqrt{1-\rho^{2}}w_{2n})^{2}+w_{3n}^{2}  &\mbox{if  } -\frac{\sqrt{1-\rho^{2}}}{\rho} w_{2n} \leq w_{1n} \leq \frac{\rho}{\sqrt{1-\rho^{2}}} w_{2n}\\
 & $  $w_{2n} \geq 0 \\
w_{3n}^{2} &\mbox{if  } w_{1n} \leq 0, w_{2n} \leq -\frac{\rho}{\sqrt{1-\rho^{2}}} w_{1n} 
\end{cases}
\]
Under $H_{0}$, the asymptotic distribution of $T_{n}$  follows $(\frac{1}{2} - \frac{\beta}{2 \pi}) \chi^{2}_{1} + \frac{1}{2} \chi^{2}_{2} + \frac{\beta}{2 \pi} \chi_{3}^{2}$ with $\beta = \cos ^{-1} (\rho)$. 
\end{theorem}

\begin{corollary}\label{corollary4}
If $S_{\xi_{1}\xi_{2},n} = 0$, then the log likelihood ratio of the sequence of experiment reduces to 
\[
\log \left ( \frac{dP_{n, \xi_{n}, \theta_{n}}}{dP_{n, \xi_{0}, \theta}} \right ) = t^\top S_{n} - \frac{1}{2} t ^\top J t + o_{p}(1),
\]
with $S_{n} = (S_{\xi_{1}^{2},n}, S_{\xi_{2}^{2},n}, S_{\theta, n}^\top)^\top \sim \mathcal{N}(0, J)$. Then the limit experiment of $\upsilon_{n}$ is based on observing one sample from $Y= t + v$ with $t = (\delta_{1}^{2}, \delta_{2}^{2}, \delta_{3}^\top)^\top \in \mathbb{R}_{+}^{2} \times \mathbb{R}^{p}$ and $v \sim \mathcal{N}(0, J^{-1})$. Proceed as in Theorem \ref{theorem4} with $u_{n} = (\tilde S_{\xi_{1}^{2},n}, \tilde S_{\xi_{2}^{2},n})^\top$ and find the corresponding Cholesky decomposition $\Lambda$ for $\Sigma_{(11)}$ and $w_{n} = \Lambda^{-1} u_{n}$. Under $H_{0}$, the asymptotic distribution of $T_{n}$ follows $(\frac{1}{2} - \frac{\beta}{2 \pi}) \chi^{2}_{0} + \frac{1}{2} \chi^{1}_{2} + \frac{\beta}{2 \pi} \chi_{2}^{2}$ with $\beta = \cos ^{-1} (\rho)$. 

\end{corollary}

\begin{remark}
When dimension gets higher, the construction of the $C(\alpha)$ test follows the similar idea. We first find residual score $\tilde S_{\xi,n}$ for $(\xi_{1}, \dots, \xi_{q})$ by projecting away the effect of the score of $\theta$. LeCam's third lemma implies that asymptotically $\tilde S_{\xi,n}$ follows $\mathcal{N}(0,\Sigma)$ under $H_{0}$ and $\mathcal{N}(\Sigma (\delta_{1}^{2}, \dots, \delta_{q}^{2}, (2\delta_{j}\delta_{k})_{j \neq k })^\top, \Sigma)$ under local alternative. The construction of the $C(\alpha)$ test is to find 
\begin{equation} \label{Tn}
T_{n} = \tilde S_{\xi,n} \Sigma^{-1} \tilde S_{\xi,n} - \inf _{\mu \in \mathcal{C}}(\Sigma^{-1} \tilde S_{\xi,n} - \mu) ^\top \Sigma (\Sigma^{-1} \tilde S_{\xi,n} - \mu) 
\end{equation}
where the cone $\mathcal{C} = \mathbb{R}_{+}^{q} \times \mathbb{R}^{q(q-1)/2} $, the space of the vector $(\delta_{1}^{2}, \dots, \delta_{q}^{2}, (2\delta_{j}\delta_{k})_{j \neq k })^\top$. We observe that $T_{n}$ is the LR statistics treating $\Sigma^{-1} \tilde S_{\xi,n}$ as the single observation in the limit experiment (See a similar idea in \citeasnoun{Silvapulle}). The $w_{n}$ worked out in Theorem \ref{theorem4} and Corollary \ref{corollary4} are explicit solution for (\ref{Tn}) when $ q=2$. For $q >2$, the solution for $\mu$, and therefore the test statistics $T_n$, can be easily found by using the R package {\tt quadprog}, \citeasnoun{quadprog}. \\

The test statistic, when dimension grows, continues to follow a mixture of $\chi^2$ distribution asymptotically under the null, albeit with more complex weights. In the simplest case, if both $J$ and $\Omega$ happen to be diagonal matrices, then all off-diagonal terms in the second-order score matrix drop and the weights take a very simple form. For $\xi \in \Xi \subset \mathbb{R}^{q}$ and let the residual score for $\xi$ be $\tilde S_{\xi, n}$ with its covariance matrix as $\Sigma_{q}$. The diagonality of $J$ implies diagonality of $\Sigma_{q}$. The optimal test statistic for $H_{0}: \xi_{1} = \dots = \xi_{q} =0$ against $H_{a}: \xi_{i} \neq 0$ for at least one $i$ is 
\begin{displaymath}
T_{n} = (0 \vee \tilde S_{\xi, n}) ^\top \Sigma_{q}^{-1} ( 0 \vee \tilde S_{\xi, n})
\end{displaymath}
Under $H_{0}$, $T_{n} \sim \sum_{i=0}^{q} \binom{q}{i} 2^{-q} \chi_{i}^{2}$. As $q$ becomes large, paying attention to the one-sided nature of the test achieves much better power performance than simply using the inner product of $\tilde S_{\xi, n}$ and the $\chi^2$ asymptotics, because the latter wastes $1-(1/2)^q$ portion of the type-I error. This point is also stressed in \citeasnoun{AkhHallin} on the optimal detection of random coefficient in autoregressive models. 
\end{remark}

\section{Examples}
In this section, we describe four examples of using the $C(\alpha)$ test for unobserved parameter heterogeneity in various models. The first Poisson regression example leads to similar test statistics already familiar in the literature. This is to illustrate that the $C(\alpha)$ test serves as a unification of many tests already available. As another example not fully elaborated here, \citeasnoun{Kiefer84} and \citeasnoun{Lancaster85} develop tests for parametric heterogeneity in Cox proportional hazard model both of which can be formulated as $C(\alpha)$ tests. Some of these familiar tests are derived under very specific assumptions on the heterogeneity distribution $F$. As we have already noted, this is not necessary as long as some very mild moment conditions are satisfied. All the other three examples are multi-dimensional cases, as this is the area where we think the limit experiment and the $C(\alpha)$ test offers most interesting departures from existing work. \\


\subsection{Tests for overdispersion in Poisson Regression} \label{poicase}
Overdispersion tests for Poisson models constitute the most common example on test of parameter heterogeneity. Such a test was proposed in \citeasnoun{Fisher50} and also serves as the motivating example in \citeasnoun{NeymanScott}. We will consider  two distinct versions of the test for unobserved heterogeneity in the conditional mean function of the Poisson regression model. \\

\subsubsection{Second Moment Test} \label{secondmoment}
Suppose we have $(Y_{1}, \dots, Y_{n})$ as i.i.d. random variables follow Poisson distribution with mean parameter $\lambda_{i}$. We further assume that
\[
\lambda_{i} = \lambda_{0i} e^{\xi U_{i}} = \exp(x_{i}^{\prime}\beta + \xi U_{i})
\]
where $U_{i}$ are i.i.d. with distribution $F$, zero mean and unit variance. We have set $\tau$ to be 1 without loss of generality. The $x_{i}$'s are covariates of the Poisson regression model including an intercept term. These covariates could be viewed as observed heterogeneity in the mean function, while $U_{i}$, since it is not explained by the covariates, is unobserved heterogeneity. Thus, the intercept coefficient, $\beta_{0}$, given the assumed form for $\lambda_i$, can be regarded as a random coefficient. We would like to test $H_{0}: \xi=0$ against $H_{a}: \xi \neq 0$ with $\beta$ as the unspecified nuisance parameters. Since the first-order score with respect to $\xi$ vanishes, this problem falls into the framework we considered in Section 2. \\

With some straightforward calculation and the nuisance parameters replaced by their MLEs, we find the $C(\alpha)$ test statistic as 
\[
Z_{n} = \frac{\sum_{i} [(y_{i} - \exp(x_{i}^{\prime} \hat \beta))^{2} - \exp(x_{i}^{\prime} \hat \beta)]}{\sqrt{2 \sum_{i} \exp (2 x_{i}^{\prime} \hat \beta)}}
\]


We call this a second moment test because $Z_{n}$ is essentially comparing the sample second moment with the second moment for the Poisson model under $H_{0}$. We reject $H_{0}$ when $(0 \vee Z_{n}) ^{2} > c_{\alpha}$ with $c_{\alpha}$ as the critical value from the mixture of $\chi^{2}$. \\

\begin{remark}
The $C(\alpha)$ test constructed above is identical to the first test statistic proposed in \citeasnoun{Lee86} for overdispersion in Poisson regression models. In his derivation, Lee assumed that the Poisson mean parameter, $\lambda_{i}$, follows a Gamma distribution with certain mean-variance ratio. The Poisson-Gamma compound distribution then leads to a negative binomial model. As Lee noted (p.700), the same test statistic can also be derived under some other distribution in addition to the Gamma distribution (See also \citeasnoun{DeanLawless89}). From the $C(\alpha)$ perspective, the test statistic does not depend on the distribution of $U$, as long as the moment conditions are satisfied. However, the form of the test statistic does depend on the particular specification on $\lambda_{i}$ as a function of the observed covariates and the unobservable $U_{i}$. This leads us to the next example. \\
\end{remark}

\subsubsection{Second Factorial Moment Test} \label{secondfac}
If instead, under the same setup as we have in \ref{secondmoment}, we assume, 
\[
\lambda_{i} = \lambda_{0i} \left (1+ \xi U_{i} / \sqrt{\lambda_{0i}} \right)
\]
The residual score for $\xi$ is now found to be, with $\lambda_{0i}=\exp(x_{i} ^{\prime} \beta)$,  
\[
g(y_i, \beta) = \left [y_{i}(y_{i}-1) - 2\lambda_{0i}(y_{i} - \lambda_{0i})-\lambda_{0i}^{2} \right ]/\lambda_{0i}
\]
and $\mathbb{V}(g(Y_i, \beta))=2$. Replacing $\beta$ by its restricted MLE $\hat \beta$, the locally optimal $C(\alpha)$ test is 
\[
Z_{n} = \frac{1}{\sqrt{2n}} \sum_{i} \left [y_{i}(y_{i}-1) - \hat \lambda_{0i}^{2} \right ] /\hat \lambda_{0i}
\]

The test statistic $Z_{n}$ is comparing the second sample factorial moment with that induced by the Poisson model under the null. Note that this test reduces to the second moment test if there are no covariates. Noticing again that only overdispersion is possible when deviating from the null, one-sided alternatives and the mixture of $\chi^{2}$ asymptotics is employed. \\

\subsection{Joint test for slope heterogeneity in linear regression model} \label{BP}
We consider a linear cross sectional model,  
\[
y_i = x_i ^\top \beta_i + u_i, 
\]
where $\beta_i$ is a $p \times 1$ vector and $u_i \sim IID \mathcal{N}(0, \sigma^2)$. In addition, we assume $\beta_{ki} = \beta_{k0} + \xi_k  U_{ki}$ for $k = 2, \dots, p$. Without loss of generality, we impose $U_{ki} = U_{i}$ for all $k$ and $U_i$ has mean zero and unit variance. As discussed earlier, this implies we need to include all distinct elements in the second order score matrix. Replacing nuisance parameters by their MLEs, it is easy to find the respective score for $\xi$ and for the nuisance parameters $\theta = (\beta^\top, \sigma^2)^\top$: 
\[
\begin{array}{ll}
S_{\xi, i} &=(\hat u_{i}^{2}/\hat \sigma^2 - 1) z_i /\hat \sigma^2\\
S_{\sigma^2,i} &= (\hat u_{i}^2/\hat \sigma^2-1)/2\hat \sigma^2\\
S_{\beta, i} & = \frac{\hat u_i}{\hat \sigma^2} x_i
\end{array}
\]
where $\hat u_i = y_i - x_i^\top \hat \beta$ and $z_i$ is the vector of length $p(p-1)/2$ that consists distinct elements of $x_i  x_i ^\top$. The same testing problem is considered in the seminal paper by \citeasnoun{Breusch79} who propose the LM test taking the form
\[
LM = \frac{1}{2} (\sum_i z_i f_i)^\top (\sum_i z_i z_i^\top)^{-1} (\sum_i z_i f_i)
\]
with $f_i = \hat u_{i}^{2}/\hat \sigma^2 - 1$. Under $H_{0}$, the LM statistic follows $\chi^2_{p(p-1)/2}$ asymptotically. \\

The $C(\alpha)$ test takes the same score function for $\xi$ and $\theta$, but pays explicit attention to the positivity constraints on those terms in $S_{\xi,i}$ that are inherited from the diagonal terms of $x_i x_i^\top$. We can easily find the residual score for $\xi$ as 
\[
\tilde S_{\xi,n} = \frac{1}{\sqrt{n}} \sum_i (z_i - \bar z)(\hat u_{i}^{2}/\hat \sigma^2 - 1) / \hat \sigma^2
\]
and the associated Information matrix as $\Sigma = 2 (\sum_i (z_i - \bar z)(z_i - \bar z)^\top)/N\hat \sigma^4$.  Partition $\tilde S_{\xi,n}$ and $\Sigma$ such that $\tilde S_{(1)}$ and $\Sigma_{(11)}$ correspond to the elements inherited from the diagonal elements of $x_i x_i^\top$ and proceed as in Theorem \ref{theorem4}. In the simulation section we give a comparison between the $C(\alpha)$ test and the LM test which provides some encouraging evidence of power improvement.

\subsection{Joint test for location and scale heterogeneity in Gaussian panel data model} \label{3.4}
In this example, we consider a two dimensional $C(\alpha)$ test for parameter heterogeneity in a Gaussian panel data model. The model is assumed to be 
\[
y_{it} = \mu_{i} + \sigma_{i} \epsilon_{it}
\]
with $\epsilon_{it} \sim IID \mathcal{N}(0,1)$,  $\mu_{i} = \mu_{0} + \xi_{1} U_{1i}$ and $\sigma_{i}^{2} = \sigma_{0}^{2} \exp (\xi_{2} U_{2i}) \geq 0$. For convenience, we assume the random variables $U_{ki}$ are i.i.d. with distribution $F_{k}$ for $k=1,2$. Both $U_{1}$ and $U_{2}$ have zero mean and unit variance and are assumed to be independent for simplicity.\\

The unconditional density of observing $(y_{i1}, \dots, y_{iT})$ is
\[
f_{i} = \int \int \Big (\frac{1}{2\pi \sigma_{0}^{2} \exp (\xi_{2} u_{2i})} \Big)^{T/2} \exp \left (-\sum_{t=1}^{T} \frac{(y_{it}-\mu_{0}-\xi_{1}u_{1i})^{2}}{2 \sigma_{0}^{2}\exp(\xi_{2}u_{2i})} \right)dF_{1}(u_{1i})dF_{2}(u_{2i})
\]
The respective score for $(\xi_{1}, \xi_{2})$ and the nuisance parameters $(\mu_{0}, \sigma_{0}^{2})$ are 
\[
\begin{array}{ll}
v_{1i}=\nabla_{\xi_{1}}^{2} \log f_{i} |_{\xi_{1}=\xi_{2}=0}&= (\frac{\bar y_{i.} - \mu_{0}}{\sigma_{0}^{2}/T})^{2}-\frac{1}{\sigma_{0}^{2}/T}\\
v_{2i}=\nabla_{\xi_{2}}^{2} \log f_{i} |_{\xi_{1}=\xi_{2}=0} & = (Z_{i}-\frac{T}{2})^{2}-Z_{i}\\
v_{3i}=\nabla_{\mu_{0}} \log f_{i} |_{\xi_{1}=\xi_{2}=0} & = \frac{\bar y_{i.}-\mu_{0}}{\sigma_{0}^{2}/T}\\
v_{4i} = \nabla_{\sigma_{0}^{2}} \log f_{i} |_{\xi_{1}=\xi_{2}=0} & = (Z_{i}-\frac{T}{2})/\sigma_{0}^{2}
\end{array}
\]
where $\bar y_{i.}$ is the sample mean defined as $\sum_{t=1}^{T} y_{it}/T$ and $2Z_{i}=\sum_{t=1}^{T}(y_{it}-\mu_{0})^{2}/\sigma_{0}^{2} \sim \chi^{2}_{T}$. \\

Replacing the nuisance parameters by their MLEs, the optimal $C(\alpha)$ test for $H_{0}: \xi_{1}=\xi_{2}=0$ against $H_{a}: \xi_{i} \neq 0$ for at least one $i$ is:
\[
T_{n} = (0 \vee t_{1n})^{2} + (0 \vee t_{2n})^{2}
\]
with
\[
\begin{array}{ll}
t_{1n} &= (2NT(T-1)/\hat \sigma_{0}^{4})^{-1/2} \left ( \sum_{i}(\frac{\bar y_{i.}-\hat \mu_{0}}{\hat \sigma_{0}^{2}/T})^{2}-\frac{NT}{\hat \sigma_{0}^{2}} \right )\\
t_{2n} &= (NT(T/2+1))^{-1/2} \left ( \sum_{i} (Z_{i} - T/2)^{2}-\frac{NT}{2} \right)
\end{array}
\]
We reject $H_{0}$ for $T_{n} > c_{\alpha}$ where $c_{\alpha}$ is the $(1-\alpha)$-quantile of $\frac{1}{4}\chi_{0}^{2} + \frac{1}{2} \chi_{1}^{2} + \frac{1}{4}\chi_{2}^{2}$. \\

\begin{remark}
The first component $t_{1n}$ of the test statistics may be recognized again as the test for individual effect in Gaussian panel data model proposed by \citeasnoun{Breusch80}. The second component $t_{2n}$ is equivalent to a single parameter $C(\alpha)$ test for a Gamma model with heterogenous scale parameter. (Analytical derivation details appear in the Appendix B.) The factorization provided by the Gaussian model leads to simple asymptotics of the test statistics. Introducing dependence between the random effects $U_1$ and $U_2$ will add an extra score function which is the cross term in the second order score matrix, $\nabla^{2}_{\xi_1 \xi_2} \log f_i$. In this case, we proceed as in Theorem \ref{theorem4}. Notice the above test is valid for the large $N$ fixed $T$ setting, and the local alternative for $\xi_n$ is of order $N^{-1/4}$. If $T$ also tends to infinity, then the local alternative for $\xi_n$ is of order $N^{-1/4}T^{-1/2}$. 
\end{remark}

\subsection{Test for slope heterogeneity in large panels}
Example \ref{3.4} above tests for randomness in individual location and variances. Perhaps a more realistic application is to allow for individual effects and the group-wise heteroscedasticity in the error but test for randomness in the slope coefficients. This problem has been considered in \citeasnoun{Swamy} and is recently revived in \citeasnoun{PY08} (hereafter PY). The PY test is a standardized version of \citeasnoun{Swamy} under large N large T setting. The model is assumed to be, 
\[
y_{it} = \alpha_i + x_{it} ^\top \beta_i + \epsilon_{it},
\]
with $\beta_i$ being a $p \times 1$ vector. The null hypothesis of interest is $H_{0}: \beta_i = \beta$ for all $i$ against $H_{1}: \beta_{i} \neq \beta_{j}$ for at least one pair of $i \neq j$. The PY test is 
\[
\tilde \Delta^{PY} = \sqrt{\frac{N(T+1)}{T-k-1}} \Big (\frac{N^{-1} \tilde S - k}{\sqrt{2k}} \Big)
\]
with $M_{\tau}$ being the familiar demean matrix and $\tilde S = \sum_i (\hat \beta_i - \hat \beta_{WFE})^\top X_{i}^\top M_{\tau} X_i (\hat \beta_i - \hat \beta_{WFE})/\tilde \sigma_{i}^{2}$ where $\hat \beta_i$ is the within estimator for each individual regression and $\hat \beta_{WFE}$ is the proper pooled estimator that accounts for individual specific variance $\tilde \sigma_{i}^{2}$. \citeasnoun{Su2013} gives an LM test interpretation for $\tilde S$ in the PY test that, 
\[
\tilde S = \sum_{i} \hat \epsilon_{i}^\top M_{\tau} X_{i} (X_{i}^\top M_{\tau} X_{i})^{-1}X_{i}^\top M_{\tau} \hat \epsilon_{i}/\tilde \sigma_{i}^{2}
\]
with $\hat \epsilon_{it} = M_{\tau}(y_{it} - x_{it}^\top \hat \beta_{WFE})$. This is the LM test statistic for considering the regression $\hat \epsilon_{it} = \alpha_{i} + (x_{it}- \bar x_{i})^\top \phi_{i} + \eta_{it}$ and test for $\phi_i = 0 $ for all $i$. As both $N$ and $T$ goes to infinity, with proper re-centering and standardization, the resulting PY test has a standard normal asymptotics under $H_{0}$ and the authors recommend a two-sided test for inference. \\

In the large N large T setting, we can also construct the $C(\alpha)$ score test for heterogeneity in coefficients. Assuming again $\beta_{ki} = \beta_{k0} + \xi_{k} U_{ki}$ for $k = 1, \dots, p$. The score function for $\xi$, $S_{\xi, n}$, is the distinct $p(p+1)/2$ elements of the second-order score matrix, which takes the form $\frac{1}{\sqrt{N}} \sum_{i} (X_{i} ^\top M_{\tau} \hat \epsilon_{i} \hat \epsilon_{i}^\top M_{\tau} X_{i}/\hat \sigma_{i}^{4} - X_{i}^\top M_{\tau} X_{i}/\hat \sigma_{i}^{2})$ with nuisance parameters replaced by MLEs. The elements of $S_{\xi,n}$ are asymptotically jointly normal with mean zero and covariance matrix $\Sigma$ under $H_{0}$ and by LeCam's third lemma, they jointly follow $\mathcal{N}(\Sigma t, \Sigma)$ under the local alternative ($\xi_{j,n} = \xi_j + \delta_j N^{-1/4}T^{-1/2}$, j = 1, \dots, p) with $t = (\delta_{1}^{2}, \dots, \delta_{p}^{2}, (2\delta_{j}\delta_{k})_{j\neq k})^\top$ as discussed in Section \ref{2.3.3}. Not surprisingly, given the connection to the score test shown by \citeasnoun{Chesher84}, this shares considerable similarity to the \citeasnoun{White82} Information Matrix test. However, the IM test rejects $H_{0}$ if $S_{\xi,n} \Sigma^{-1} S_{\xi,n} $ exceeds the critical value from $\chi^2_{p(p+1)/2}$ at nominal level $\alpha$, while the $C(\alpha)$ test modifies the IM test by adjusting for positivity constraints in $t$ for the respective elements in the score function. We do not repeat the steps here in applying Theorem \ref{theorem4}. In the simulation section, we compare the $C(\alpha)$, the IM test and the PY test and the results show that the $C(\alpha)$ test enjoys a power gain compared to the other two tests. It is also worth mentioning that the advantage of the $C(\alpha)$ test is that we only need to estimate under the null model. In addition, the test can be derived in the same way for large N and fixed T setting, except the local alternative for $\xi_n$ is then of order $N^{-1/4}$. \\

\section{Reparameterization and connection to the Information matrix Test} \label{4}
\subsection{Reparameterization}
A common strategy in prior literature to circumvent the irregularity, that the first-order score function is degenerately zero, is to reparameterize the model. In fact, this is the advice given in the original \citeasnoun{Neyman59} $C(\alpha)$ paper (Section 9, p. 225) and also in \citeasnoun[p. 117-118]{CoxHinkley}. For the heterogeneity tests considered in this paper in particular, \citeasnoun{Cox83} and \citeasnoun{Chesher84} adopt such a reparameterization by letting $\eta = \xi_{0} + (\xi-\xi_{0})^{2}$. Reconsidering the example in Section \ref{2.2}, without loss of generality, we set $\xi_{0}=0$ and have the density function as $p(x; \lambda_{0} + \tau \sqrt{\eta} U_{i})$. \citeasnoun{Cox83} tests for heterogeneity of $\lambda_{i}$ by testing $H_{0}: \eta =0$ against $H_{1}: \eta >0$. Under $H_0$, this is to test whether $\Var(\lambda) = 0$. \citeasnoun{Chesher84} takes the same model assuming $U_{i}$ follows a symmetric location-scale distribution. A more recent treatment, focusing on random individual effects in panel data models by \citeasnoun{Bennala12} also uses the same reparametrization but adopts a less stringent LeCam framework. \\

At first sight, reparameterization avoids the irregularity of having a degenerate score function. The first order derivative with respect to $\eta$, albeit an undefined $\frac{0}{0}$ function, can be evaluated by the l'H\^opital's rule. As long as $\mathbb{E}(U^{2})$ is non-zero, the score function is nonvanishing. The score function thus derived also involves the second derivative and is identical to the score function in the $C(\alpha)$ test using the original parameterization that $\lambda_{i} = \lambda_{0} + \tau \xi U_{i}$. However, the second order derivative for $\eta$ is unbounded unless we impose an additional moment condition on $U$, that $\mathbb{E}(U^{3})=0$ (See the derivation in the Appendix C). This condition is implicitly satisfied in \citeasnoun{Chesher84} because of the symmetry distribution assumption on $U$. \citeasnoun{Moran73} also employed this zero third moment condition and remarked that it was hard to rationalize. One explanation for this extra condition is that the original, more natural specification on the random parameter $\lambda_{i} = \lambda_{0} + \tau \xi U_{i}$ with $\xi \in \mathbb{R}$ is not equivalent to the reparameterization $\lambda_{i} = \lambda_{0} + \tau \sqrt{\eta} U_{i}$ with $\eta \in \mathbb{R}_{+}$ unless $U$ has a symmetric distribution. When symmetry does not hold for the distribution of $U$, the likelihood does not obtain a proper expansion around $\eta$. As we have seen, the $\xi$ parameterization has the advantage that no symmetry or higher moment conditions on $U$ are necessary. \\

\subsection{Connection to the Information Matrix test}\label{4.2}
\citeasnoun{Chesher84} was the first to point out that \possessivecite{White82} Information Matrix (IM) test is a score test for unobserved heterogeneity. Since \citeasnoun{Chesher84} can be viewed as a reparameterized $C(\alpha)$ test, it is of interest to investigate the connection between the $C(\alpha)$ test for heterogeneity in general and the IM test. We show that the $C(\alpha)$ test for heterogeneity nests the IM test as a special case. \\

Take again the example in Section \ref{2.2}, $Y_{1}, \dots Y_{n}$ are i.i.d. random variables each with density function $p(y; \lambda_{i})$. The parameter $\lambda_{i}$ is a random parameter and we assume it now takes a more general form $\lambda_{i} = \lambda_{0} + \xi k(\lambda_{0})U_{i}$ to incorporate both additive and multiplicative specifications. For example, if $k(\lambda_{0})=1$, we have the additive form $\lambda_{i} = \lambda_{0} + \xi U_{i}$, while if $k(\lambda_{0})= \lambda_{0}$, then the multiplicative form. The function $k(\lambda_{0})$ thus allows flexible specification for the random parameter. \\

For simplicity and to fix ideas, we first assume $\lambda_{0}$ is known. Theorem 1 then implies the following expansion of the log-likelihood function, provided that $\xi_{n} = O(n^{-1/4})$, 
\[
l = \sum_{i} \log \int p(y_i; \lambda_{i})dF(u) = \sum_{i} \log p (y_i; \lambda_{0}) + \frac{1}{2} \xi_{n}^{2} \mathbb{E}(U_{i}^{2}) \sum_{i} k(\lambda_{0})^{2} \frac{\nabla_{\lambda}^{2}p(y_i; \lambda_{0})}{p(y_i; \lambda_{0})} + O_{P}(1)
\]
The first order derivative of $l$ with respect to $\xi_{n}$ is zero evaluated under $\xi_n=0$, and the second-order score is
\[
\frac{\partial ^{2}}{\partial \xi_{n}^{2}}l|_{\xi_{n}=0}  = \sum_{i} k(\lambda_{0})^{2} \frac{\nabla_{\lambda}^{2}p(y_i; \lambda_{0})}{p(y_i; \lambda_{0})}.
\]
If $\lambda_{0}$ is unknown, we find the corresponding score for $\lambda_{0}$ and take the projection step to get the $C(\alpha)$ test. This is very close to the approximation in \citeasnoun{Cox83} except we allow for a more flexible variance function for the random parameter $\lambda_{i}$, as $\xi^{2}\mathbb{E}(U_{i}^{2})k(\lambda_0)^{2}$. In a regression model with covariates, $\lambda_{0}$ will then be a function of the covariates with coefficients $\beta$. \\

\possessivecite{White82} Information Matrix test under regression setting, on the other hand, is constructed based on the following moment conditions:
\[
\mathbb{E}\left [ vech \left ( \nabla_{\beta}^{2} \log p(y; \lambda_{0}(x_{i}, \beta)) + \nabla_{\beta} \log p(y; \lambda_{0}(x_{i}, \beta))\nabla_{\beta}^\top \log p(y; \lambda_{0}(x_{i}, \beta)) \right) \right ] =0
\]
where $vech$ is the operator which stacks the elements in the lower triangular part of a symmetric matrix. Using the chain rule, we see that the IM test statistic uses the following sample analogue of the moment condition
\[
IM= \sum_{i}\Big[ \frac{\nabla_{\lambda}^{2} p(y; \lambda_{0}(x_{i}, \beta))}{p(y; \lambda_{0}(x_{i}, \beta))}  \nabla_{\beta} \lambda_{0}(x_{i}, \beta) \nabla_{\beta}^\top \lambda_{0}(x_{i}, \beta) + \frac{\nabla_{\lambda}p(y; \lambda_{0}(x_{i}, \beta))}{p(y; \lambda_{0}(x_{i}, \beta))} \nabla_{\beta}^{2} \lambda_{0}(x_{i}, \beta) \Big]
\]

There are various forms for the IM test in the literature (see \citeasnoun{DavMac}), we focus on the efficient score version, in which all the nuisance parameters are replaced by their restricted MLEs. For the $C(\alpha)$ test to be equivalent to the efficient score version of the IM test, it is sufficient to have the following two identities:
\[
\begin{array}{cc}
C \nabla_{\beta} \lambda_{0}(x_{i}, \beta) \nabla_{\beta}^\top \lambda_{0}(x_{i}, \beta) = k(\lambda_{0}) k (\lambda_{0})^\top\\
\sum_{i} \frac{\nabla_{\lambda}p(y; \lambda_{0}(x_{i}, \beta))}{p(y; \lambda_{0}(x_{i}, \beta))} \nabla_{\beta}^{2} \lambda_{0}(x_{i}, \beta)=0
\end{array}
\]
where $C$ is a non-zero constant. We give several examples below as illustrations. \\

\begin{example}
\textit{Normal regression with $Y_{i} \sim \mathcal{N}(\mu_{i}, 1)$, where $\mu_{i} = \mu_{0i} + \xi k(\mu_{0i}) U_{i} $ and $\mu_{0i}=x_{i}^{\prime}\beta$.}\\

\noindent Note that $\nabla_{\beta}\mu_{0i}\nabla_{\beta}^\top\mu_{0i}=x_{i}x^\top_{i}$ and $\nabla_{\beta}^{2}\mu_{0i}=0$. Considering only the IM test based on the intercept term, it is equivalent to the $C(\alpha)$ test for heterogeneity in $\beta_{0}$ if $k(\mu_{0i}) = C \neq 0$. If considering all elements in the IM test, the equivalence holds if $x_{i} x^\top _{i} = k(\mu_{0i}) k(\mu_{0i})^\top$. In this case, the $C(\alpha)$ test is  multivariate, testing for homogeneity for all coefficients $\beta$ in $\mu_{0i}$. 
\end{example}

\begin{example}\label{poi}
\textit{Poisson regression with $Y_{i} \sim Poi(\lambda_{i})$, where $\lambda_{i} = \lambda_{0i} + \xi k(\lambda_{0i}) U_{i}$ and $\lambda_{0i}=\exp(x_{i}^{\prime}\beta)$.} \\

\noindent Considering only the IM test for the intercept term,  we have $\nabla_{\beta_{0}} \lambda_{0i} = \nabla_{\beta_{0}}^{2} \lambda_{0i} = \lambda_{0i}$. If $\beta$'s are replaced by their MLEs, the second identity for equivalence holds because the normal equation for the MLE of $\beta_{0}$ gives 
\[
\sum_{i} \frac{\nabla_{\lambda}p(y; \lambda_{0i})}{p(y; \lambda_{0i})} \nabla_{\beta_{0}}^{2} \lambda_{0i} = \sum_{i} \frac{\nabla_{\lambda}p(y; \lambda_{0i})}{p(y; \lambda_{0i})} \nabla_{\beta_{0}} \lambda_{0i} =0
\]
Therefore, the IM test is equivalent to the $C(\alpha)$ test if $k(\lambda_{0i})= \lambda_{0i}$ which is satisfied for the multiplicative alternative $\lambda_{i} = \lambda_{0i}(1+  \xi U_{i})$. This specification is a first order linear approximation of the alternative form $\lambda_i=\lambda_{0i}\exp(\xi U_{i})$ for small $\xi$, which leads to the second moment test for the Poisson regression model as discussed in Section \ref{secondmoment}. There are of course many other possible specifications for the conditional mean function of $\lambda_{0i}$ which would lead to other equivalence conditions. We do not delve into further details here, but refer the readers to \citeasnoun[Chapter 5]{CamTri} and \citeasnoun{Dean92} for more elaborated discussions on count data models. 
\end{example}

In summary, when the model contains covariates, the functional form of the $C(\alpha)$ test is equivalent to the IM test only under a particular alternative specification, provided that the nuisance parameters are also replaced by their corresponding restricted MLEs. When the model does not contain covariates, the IM test will always be equivalent to the $C(\alpha)$ test because the function $k(\lambda_{0})$ is no longer individual specific and can be factored out as a constant from the score function. It will then be cancelled when we rescale the score by its standard deviation to form the $C(\alpha)$ test statistic. To see more clearly how different specification affect the power performance of various testing procedures discussed here, especially in cases where the IM test no longer serves as an optimal test, we conduct a Monte Carlo simulation in the next section. It is also important to deviate from the common practice in using the $\chi^2$ asymptotics for the IM test or the LM test for heterogeneity. The simulation shows that overlooking the intrinsic one-sidedness of alternatives sacrifices power.

\section{Simulation}
We first revisit the Poisson regression model to illustrate the points made in Section \ref{4.2}. As discussed in Example \ref{poi} and also in Section \ref{poicase}, when $k(\lambda_{0i})$ takes different functional forms, one finds different optimal test statistics. For two different data generation processes, we compare three testing procedures: the second moment test and the second factorial moment test, both are one-sided tests and use critical value from a mixture of $\chi^2$ and the information matrix test, using critical values from the $\chi^2$ distribution. The first experiment generates data from a Poisson regression model with the conditional moment function as $\lambda_{i} = \lambda_{0i} +\tau  \xi \lambda_{0i} U_{i}$ and the second with $\lambda_{i} = \lambda_{0i} +\tau  \xi \sqrt{\lambda_{0i}} U_{i}$. In both cases $\lambda_{0i} = exp(\beta_{0} + \beta_{1} x_{i})$ and $\tau \xi U_{i}$ has a mixture distribution taking value $1.5h$ with probability $2/3$ and $-3h$ with probability $1/3$. We consider 21 distinct values of $h$ equally spaced and the design of $X$ is fixed for all experiments as a sample drawn from a standard normal distribution. Using other $X$ designs does not change the conclusions. The sample size for all power comparison is fixed at 500 with 10000 replications. \\

\begin{figure}[h!]
\centering 
\includegraphics[scale = 0.45]{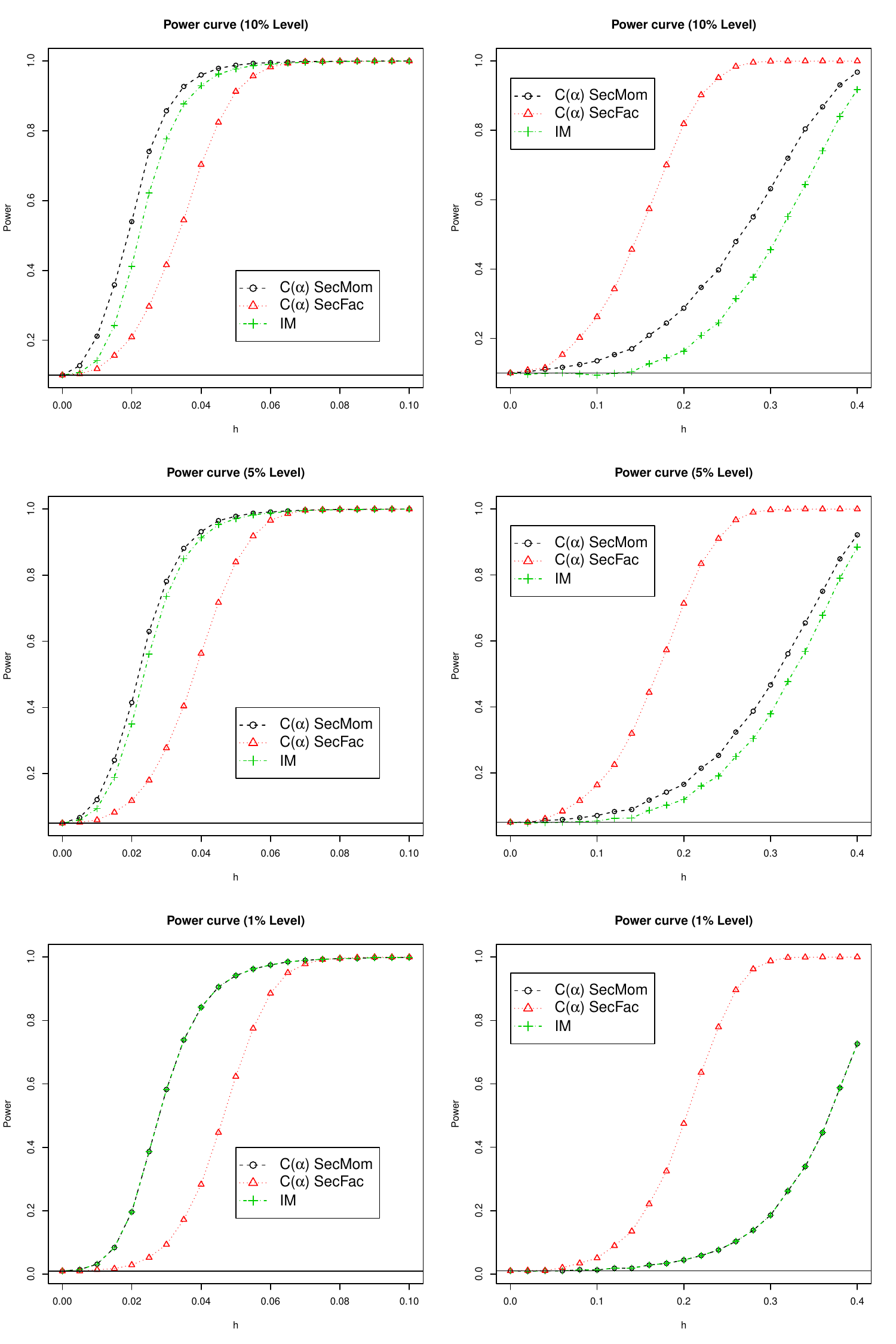}
 \caption{Power Comparison of Unobserved Heterogeneity Test for Poisson Regression Model: The left panel corresponds to the first experiment and the right panel to the second. The dotted curve corresponds to the power curve of the second moment test, the curve with triangle signs for the second factorial moment and the crossed curve for the IM test of the intercept term.}
\label{power}
\end{figure}

In the left panel of Figure \ref{power}, both the second moment test and the information matrix test performs uniformly better than the second factorial moment test. This is to be expected since the second moment test is the optimal test derived using the $C(\alpha)$ framework. The IM test using just the element for the intercept term has an identical test function as the second moment test, but using the one-sided test with mixture of $\chi^2$ critical value gives better power, especially for the $10\%$ level case. On the other hand, the second factorial moment test is superior under the second experiment, although the power for the second moment test and the IM test also converges to unity albeit much more slowly. It is documented in the literature that the IM test has poor size in small samples (\citeasnoun{ChesherSpady}) and the $C(\alpha)$ test may be subject to similar criticism. We use size-corrected critical values as suggested by \citeasnoun{Horowitz94}. \\

\begin{figure}[h!]
\centering 
\includegraphics[scale = 0.5]{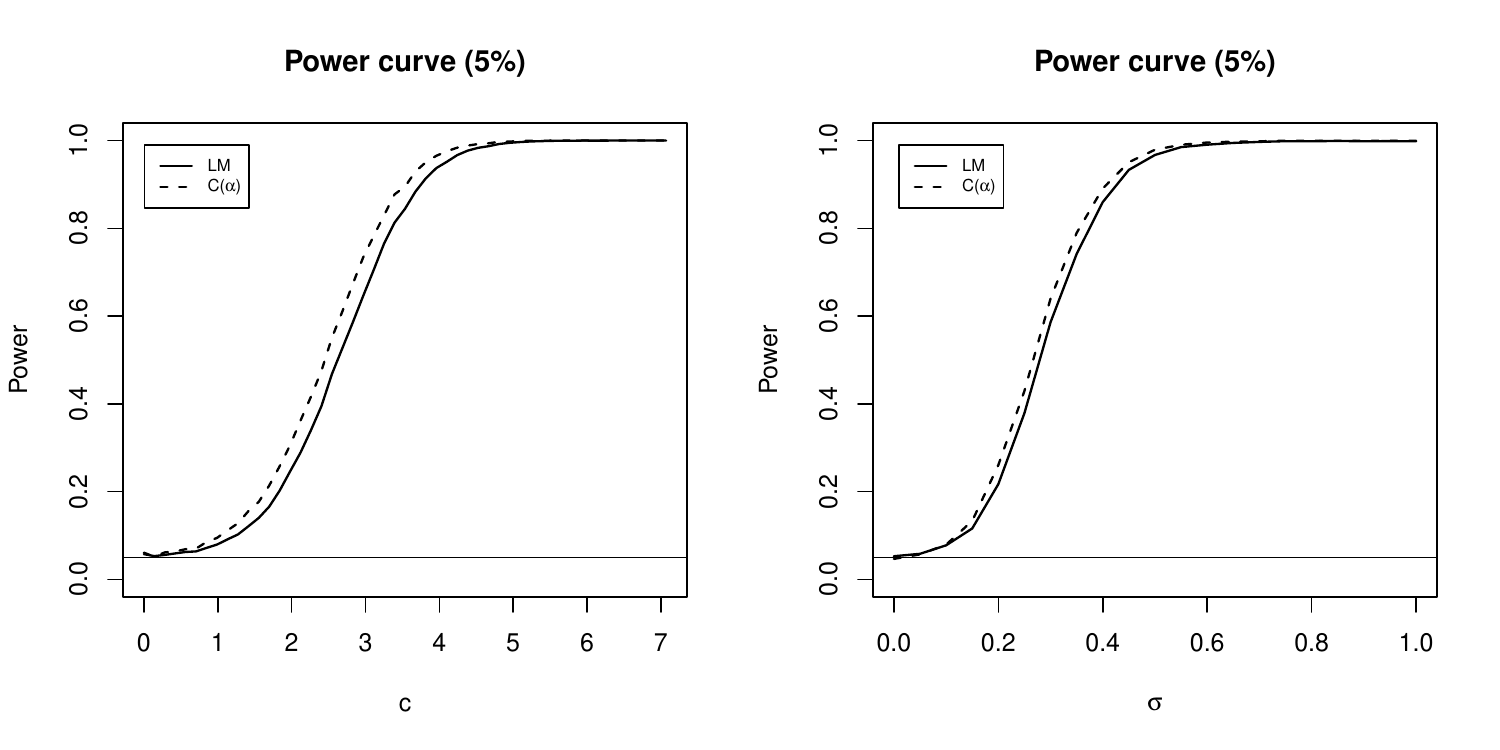}
 \caption{Power Comparison of Slope Heterogeneity Test for Linear Regression Model: The left figure corresponds to the first experiment and the right to the second. The dotted curve corresponds to the power curve of the $C(\alpha)$ test based on the mixture of $\chi^2$ asymptotics and the solid curve for the LM test based on the $\chi^{2}_{3}$ asymptotics. }
\label{LMtest}
\end{figure}

We then conduct a power comparison between the $C(\alpha)$ test and the \citeasnoun{Breusch79} LM test for Example \ref{BP}. Two different alternative $\beta_{i}$ distributions are considered. The first one assumes $\beta_{i}$ takes value $0$ for $i = 1, \dots, N/2$ and $c N^{-1/4}$ for $i = N/2+1, \dots, N$.  We let  $c$ take 51 distinct values equally spaced from $0$ to $\sqrt{50}$. The second case assumes $\beta_i \sim \mathcal{N}(0,\sigma^2)$ with $\sigma$ taking 21 distinct values from $0$ to $1$. For simplicity, we consider the case with dimension two, where both $x$ covariates are standard normal variables. The sample size is fixed at 500 with 10000 replications. Figure \ref{LMtest} presents the power curve for the 5 \% nominal level. The first experiment has a slightly bigger power gain compared to the second, but in both cases, the $C(\alpha)$ test dominate the power curve of the LM test based on the usual $\chi^2$ asymptotics. \\

Lastly, we compare the $C(\alpha)$ test with the \citeasnoun{PY08} test and the Information Matrix test for a Gaussian panel data model. As shown by \citeasnoun{PY08} Table 1, their standardized Swamy test has very nice size and power performance compared to some other existing tests, i.e. the Hausman test and the original \citeasnoun{Swamy} test for a wide range of $N$ and $T$. We consider a panel data model with two exogenous regressors and normal errors with individual variances. Table \ref{tab: sizepower} reports the size and power for the $C(\alpha)$ test, the PY test and the Information matrix test. Both the PY test and the $C(\alpha)$ test has correct size and the IM test is slightly on the conservative side. For all $N$ and $T$ combinations, the $C(\alpha)$ test has a significant power gain. \\

\begin{figure}[h!]
\centering 
\includegraphics[scale = 0.5]{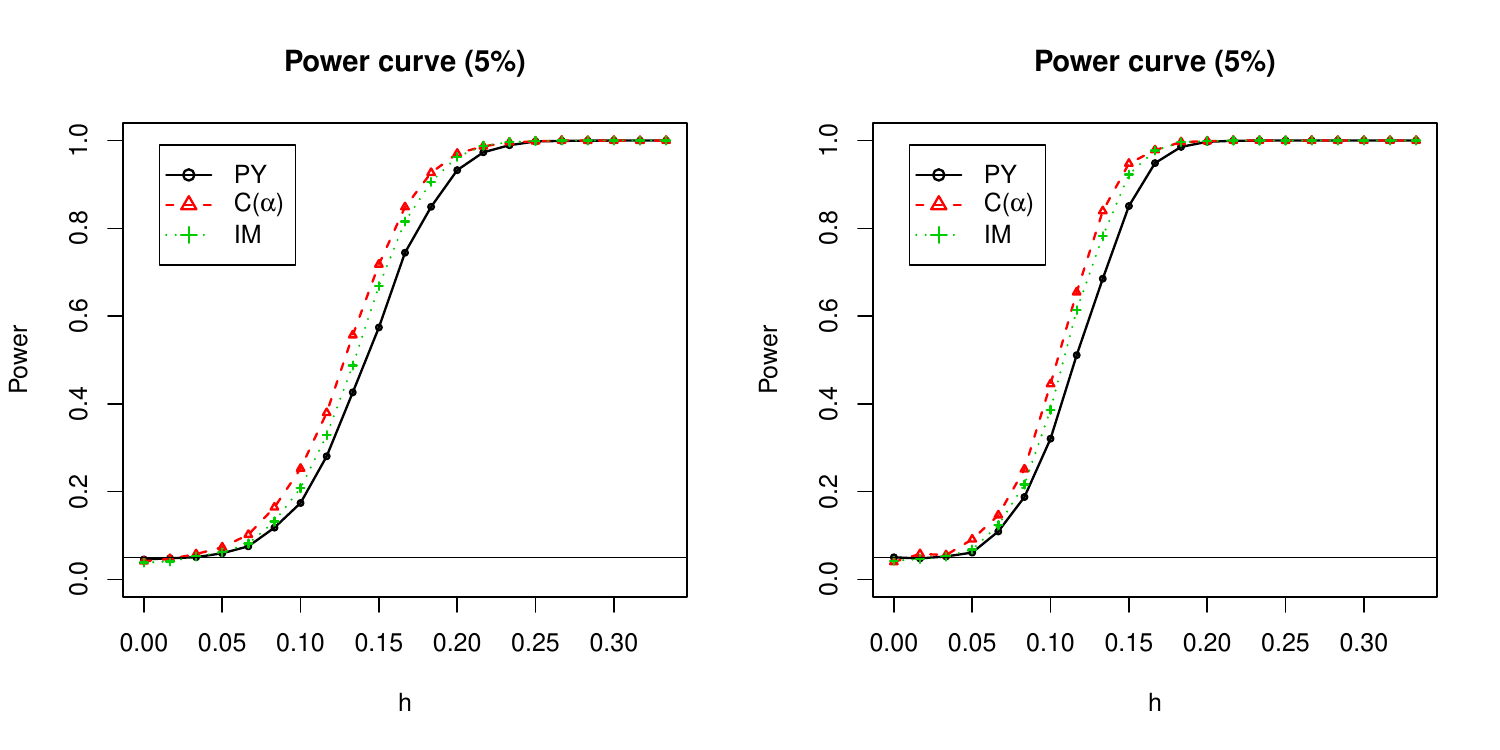}
 \caption{Power Comparison of Slope Heterogeneity Test for Gaussian Panel Data Model: The left figure corresponds to the first experiment and the right to the second for different values of $h$. Data are generated as $y_{it} = \alpha_{i} + x_{it}^\top \beta_i + \epsilon_{it}$ with $\alpha_{i} \sim U(0,1)$ and $\epsilon_{it} \sim IID \mathcal{N}(0, \sigma_{i}^{2})$ and $\sigma_{i}^{2} \sim U(1,2)$. Both regressors are normal variable with mean zero and standard deviation 0.5. The solid line with circles is the power curve for the PY test, the crossed curve for the Information Matrix test with $\chi^2$ asymptotics and the curve with triangle signs for the $C(\alpha)$ test with mixture of $\chi^2$ asymptotics. }
\label{PYtest}
\end{figure}

For a further power comparison, we consider with the same model as above for two different $\beta_i$ distributions. For 21 distinct equally spaced values of $h \in [0, 1/3]$, the first example assumes $\beta_{i} \sim \mathcal{N}(0, h^2)$ and the second assumes $\beta_{i}$ taking two possible values $\{ 1-h,1+2h\}$ with probability $(2/3, 1/3)$. Results are presented in Figure \ref{PYtest}. The sample size is fixed for each experiment at $T = 50$ and $N = 100$ with 5000 replications. The $C(\alpha)$ test again exhibits encouraging power improvement compared to the other two tests uniformly for all $h$ values. \\

\begin{table}[h]
\centering
\begin{tabular}{rrrrrrrr} 
  \hline
& &   \multicolumn{3}{c}{Size} & \multicolumn{3}{c}{Power}\\
  \cline{3-5} \cline{6-8}
  T & N & PY & $C(\alpha)$ & IM & PY & $C(\alpha)$ & IM \\ 
\hline
$ 20$&$ 30$&$0.053$&$0.040$&$0.032$&$0.064$&$0.072$&$0.055$\tabularnewline
$ 30$&$ 30$&$0.049$&$0.046$&$0.034$&$0.082$&$0.109$&$0.095$\tabularnewline
$ 50$&$ 30$&$0.046$&$0.054$&$0.041$&$0.137$&$0.215$&$0.179$\tabularnewline
$100$&$ 30$&$0.045$&$0.048$&$0.042$&$0.434$&$0.604$&$0.538$\tabularnewline
$ 20$&$ 50$&$0.045$&$0.042$&$0.037$&$0.070$&$0.070$&$0.063$\tabularnewline
$ 30$&$ 50$&$0.046$&$0.038$&$0.033$&$0.112$&$0.154$&$0.125$\tabularnewline
$ 50$&$ 50$&$0.045$&$0.052$&$0.041$&$0.263$&$0.379$&$0.325$\tabularnewline
$100$&$ 50$&$0.047$&$0.052$&$0.044$&$0.625$&$0.738$&$0.721$\tabularnewline
$ 20$&$100$&$0.046$&$0.040$&$0.034$&$0.090$&$0.099$&$0.084$\tabularnewline
$ 30$&$100$&$0.051$&$0.039$&$0.040$&$0.172$&$0.224$&$0.195$\tabularnewline
$ 50$&$100$&$0.046$&$0.045$&$0.047$&$0.406$&$0.570$&$0.484$\tabularnewline
$100$&$100$&$0.044$&$0.042$&$0.046$&$0.886$&$0.946$&$0.945$\tabularnewline
$ 20$&$200$&$0.049$&$0.035$&$0.033$&$0.151$&$0.171$&$0.140$\tabularnewline
$ 30$&$200$&$0.045$&$0.044$&$0.036$&$0.317$&$0.395$&$0.336$\tabularnewline
$ 50$&$200$&$0.045$&$0.047$&$0.041$&$0.672$&$0.800$&$0.773$\tabularnewline
$100$&$200$&$0.048$&$0.048$&$0.043$&$0.993$&$0.998$&$0.999$\tabularnewline
   \hline\\
\end{tabular}
\caption{Size and Power comparison between the PY test, the Information Matrix test and the $C(\alpha)$ test for different N and T. Data are generated as $y_{it} = \alpha_{i} + x_{it}^\top \beta_i + \epsilon_{it}$ with $\alpha_{i} \sim U(0,1)$ and $\epsilon_{it} \sim IID \mathcal{N}(0, \sigma_{i}^{2})$ and $\sigma_{i}^{2} \sim U(1,2)$. Both regressors are $\mathcal{N}(0, 0.5^2)$. Under the null, $\beta_{1i} = \beta_{2i} = 1$ for all $i$ and under the alternative, $\beta_{1i} = \beta_{2i} \sim \mathcal{N}(1, 0.15^2)$. The PY test is based on a two sized $\mathcal{N}(0,1)$ test, the IM test is based on $\chi_{3}^{2}$ test and the $C(\alpha)$ test on a mixture of $\chi^2$ test. All tests are conducted at 5\% nominal level with 5000 replications.}
\label{tab: sizepower}
\end{table}

\section{Conclusion}
We have shown that Neyman's $C(\alpha)$ test provides a unified approach to testing for neglected heterogeneity in parametric models. The irregularity encountered in these testing problems, that the score function is identically zero, can be circumvented by defining a second-order score function. Optimality of this new score function is established by formulating the problem in LeCam's LAN framework and examining the associated limit experiment. This framework provides neater regularity conditions in the irregular problem as compared to classical approach in \citeasnoun{Neyman59}. The multi-dimensional extension suggests a modification on the usual $\chi^2$ test that leads to power improvement in many applications. \\

The $C(\alpha)$ test inherits the chief merit of the score test, computation is made easy under the null model. In contrast, the likelihood ratio test, in face of the generally unknown heterogeneity distribution $F$, is computationally challenging. We have also seen that the $C(\alpha)$ test has local power against a wide class of alternatives, that allows us to avoid strict parametric assumptions on $F$, relying instead on weaker moment conditions. A further advantage of the LeCam framework is that it enables us to dispense with symmetry and higher order moment conditions that have been employed in earlier work. \\

A straightforward generalization of the theorems in Section 2 would be to incorporate density functions that allow the first $(m-1)$ logarithmic derivatives to vanish. \citeasnoun{Rotnitzky} also discuss estimation problems in this general case under classical MLE type of conditions. In such cases, we can define the $m^{th}$ order derivative of the log density as the score function and require the Pitman-type local alternative to be of order $n^{-1/2m}$. LeCam's DQM condition needs to be modified by raising the corresponding elements in the expansion to $m^{th}$ power, as we did for $m=2$ in Definition 1. It is curious to observe that only when $m$ is an even integer is the test required to be one-sided and reparameterization is not advisable. When $m$ is odd, we can use reparameterization to transform the irregular problem back to a regular case, without imposing additional restrictions (i.e. symmetry of the distribution $F$). \\

A drawback of the $C(\alpha)$ test, as reflected in \citeasnoun{Neyman79}, is that asymptotic optimality of the test is only established under local alternatives. The approximation of the power function, which is characterized by the asymptotic behavior of the test statistics under such alternatives, relies on $n$ tending to infinity and the parameter $\xi_{n}$ converging to the null value $\xi_{0}$. The behavior of the power function for finite samples or fixed alternatives is largely unknown. Some finite sample correction like those pursued in \citeasnoun{Honda88} and \citeasnoun{ChesherSpady} is left for future work. 

\bibliography{Calpha}

\begin{appendix}
\section{Proof of theorems}
\noindent Before proceeding to the proof for \textbf{Theorem 1}, we first prove the following lemma as an adaption to \citeasnoun[Lemma 1]{Pollard97}. Denote $f_{n}=\sqrt{p(x_{i}; \xi_{n}, \theta_{n})}$ and $f_{0}= \sqrt{p(x_{i}; \xi_{0}, \theta)}$. Let $v_{\xi}$ and $v_{\theta}$ be shorthand for $v_{\xi}(x_{i})$ and $v_{\theta}(x_{i})$. Let $\| \cdot \|$ be $\mathcal{L}_{2}(\mu)$-norm and $\left \langle \cdot , \cdot \right \rangle$ be the inner product. If it contains a vector, then it is defined as the vector of inner product for each elements. Further, let $r_{n}(x_{i}, \xi_{n}, \theta_{n})= f_{n}-f_{0}- h_{n}^\top v(x_{i})$ and denote $R_{i} = r_{n}(x_{i}, \xi_{n}, \theta_{n})/f_{0}$. 
\begin{lemma}
Under Assumption 1 and the modified DQM condition, we have the following:
\begin{enumerate}
\item  $\sum_{i}R_{i}^{2} = o_{p}(1)$ \label{eq:s1}
\item $\mathbb{E}(v(X)/f_{0})=0$ \label{eq:mean0}
\item $2 \sum_{i} R_{i} = -\frac{1}{4} t ^\top J t + o_{P}(1)$ \label{eq:s2}

\item $n^{-1/2}\sum_{i} R_{i} v_{\xi}/f_{0}= o_{P}(1)$, $n^{-1/2}\sum_{i} R_{i}  v_{\theta}/f_{0} = o_{P}(1)$ \label{eq:s3}
\item $\underset{1 \leq i \leq n}{\max} |R_{i}| = o_{P}(1)$ \label{eq:s4}
\item $\underset{1 \leq i \leq n}{\max} |\frac{2}{\sqrt{n}} \frac{v_{\xi}}{f_{0}}| = o_{P}(1)$, $\underset{1 \leq i \leq n}{\max} | \frac{2}{\sqrt{n}} \frac{v_{\theta}}{f_{0}}|= o_{P}(1)$ \label{eq:s5}\\

\end{enumerate}
\end{lemma}

\textbf{Proof of (\ref{eq:s1})}.  
Under the modified DQM condition, the Markov inequality yields,
\[
\begin{array}{ll}
\mathbb{P}(\sum_{i} R_{i}^{2} > \epsilon) & \leq \epsilon^{-2} n \mathbb{E}(R_{1}^{2})\\
& = \epsilon^{-2} n \int r_{n}^{2}(x; \xi_{n}, \theta_{n}) d\mu(x) \to 0.
\end{array}
\]

\textbf{Proof of (\ref{eq:mean0}) and (\ref{eq:s2})}.  
Since both $f_{n}$ and $f_{0}$ are objects with $\mathcal{L}_{2}(\mu)$-norm 1 
\[
\begin{array}{ll}
0 & = || f_{n} ||^{2}_{\mu, 2} - ||f_{0}||^{2}_{\mu, 2} \\
&=(\xi_{n}-\xi_{0})^{4} ||v_{\xi}||^{2}_{\mu, 2} + (\theta_{n}-\theta)^\top ||v_{\theta}||^{2}_{\mu, 2} (\theta_{n}-\theta) + ||r_{n}||^{2}_{\mu, 2} + 2 \left \langle (\theta_{n}-\theta)^\top v_{\theta} , r_{n} \right \rangle\\
& + 2 (\xi_{n}-\xi_{0})^{2} (\theta_{n}-\theta) ^\top \left \langle v_{\theta} , v_{\xi} \right \rangle  + 2 (\xi_{n}-\xi_{0})^{2} \left \langle v_{\xi} , r_{n} \right \rangle + 2(\xi_{n}-\xi_{0})^{2} \left \langle f_{0} , v_{\xi} \right \rangle \\
& + 2(\theta_{n}-\theta)^\top \left \langle f_{0} , v_{\theta} \right \rangle + 2 \left \langle f_{0} , r_{n} \right \rangle
\end{array}
\]
Let $\{ \theta_{n}, \xi_{n} \}$ be sequences such that $\theta_{n} - \theta = O(n^{-1/2})$ and $(\xi_{n} - \xi_{0})^{2} = O(n^{-1/2})$. Note that by Cauchy-Schwarz inequality and the fact that both $v_{\xi}$ and $v_{\theta}$ are square integrable with respect to measure $\mu$ by assumption, $\left \langle v_{\xi} , r_{n} \right \rangle=o(1/\sqrt{n})$ and $\left \langle v_{\theta} , r_{n} \right \rangle = o(1/\sqrt{n})$. Therefore, the third, fourth and the sixth terms are of order $o(1/n)$. The first, second and fifth terms are of order $O(1/n)$. The ninth term is of order $o(n^{-1/2})$ by Cauchy-Schwarz inequality. The seventh and eighth term are both of order $O(1/\sqrt{n})$, but in order for the identity to hold, they must be of smaller order to balance with other terms. For this to happen, we must have 
\[
\left \langle f_{0} , v_{\xi} \right \rangle = \left \langle f_{0} , v_{\theta} \right \rangle =0
\]
This proves (\ref{eq:mean0}) since $0=\left \langle f_{0} , v_{\xi} \right \rangle=\mathbb{E}(v_{\xi}(X)/f_{0})$. Similar argument shows $\mathbb{E}(v_{\theta}(X)/f_{0})=0$. Hence, 
\[
\begin{array}{ll}
2 \left \langle f_{0} , r_{n} \right \rangle &= -(\xi_{n}-\xi_{0})^{4} ||v_{\xi}||^{2}_{\mu, 2} - (\theta_{n}-\theta)^\top ||v_{\theta}||^{2}_{\mu, 2}(\theta_{n}-\theta) \\
&- 2(\xi_{n}-\xi_{0})^{2}(\theta_{n}-\theta)^\top \left \langle v_{\xi} , v_{\theta} \right \rangle + o(1/n) \\
& = -\frac{1}{4n} t^\top J t + o(1/n)
\end{array}
\]
with $t^\top = (\delta_{1}^{2}, \delta_{2}^\top)$. \\

Since $\mathbb{V}(2 \sum_{i}R_{i})$ is bounded above by $4\sum_{i} \mathbb{E}(R_{i}^{2})$, which goes to 0 from (\ref{eq:s1}), we have 
\[
\begin{array}{ll}
2\sum_{i}R_{i} & = 2n \mathbb{E}(R_{1}) + o_{P}(1)\\
& = 2 n \left \langle f_{0} , r_{n} \right \rangle + o_{P}(1)\\
&= 2n \left (-\frac{1}{8n} t^\top J t + o(1/n) \right) + o_{P}(1)\\
&= -\frac{1}{4} t^\top J t + o_{P}(1)
\end{array}
\]

\textbf{Proof of (\ref{eq:s3})}.  
By H\"older's inequality, 
\[
\sum_{i}R_{i} \frac{2}{\sqrt{n}} \frac{v_{\xi}}{f_{0}}  \leq \sqrt{\sum_{i}R_{i}^{2} \sum_{i} (\frac{2}{\sqrt{n}} \frac{v_{\xi}}{f_{0}})^{2}} = o_{P}(1)O_{P}(1) = o_{P}(1)
\]
Similar argument admits the second result. 

\textbf{Proof of (\ref{eq:s4})}.  
\[
\mathbb{P}(\underset{1 \leq i \leq n}{\max} |R_{i}| > \epsilon)  \leq n \mathbb{P}(|R_{i}|^{2} > \epsilon^{2})  \leq \epsilon^{-2} n \mathbb{E}(R_{i}^{2}) \to 0
\]

\textbf{Proof of (\ref{eq:s5})}.  
\[
\begin{array}{ll}
\mathbb{P}(\underset{1 \leq i \leq n}{\max} |2 v_{\xi}/f_{0}| > \epsilon \sqrt{n}) & \leq n \mathbb{P}(|2 v_{\xi} /f_{0}| > \epsilon \sqrt{n}) \\
& \leq \epsilon^{-2} \mathbb{E}\left ((2 v_{\xi}(X_{1})/f_{0})^{2} \right) \mathbb{I}_{[|2 v_{\xi}/f_{0}| > \epsilon \sqrt{n}]} \to 0
\end{array}
\]
Similar argument admits the second statement. \hfill{\rule{2mm}{2mm}} \\

\textbf{Proof of Theorem \ref{theorem1}} 
We consider $\xi_{n} = \xi_{0} + \delta_{1} n^{-1/4}$ and $\theta_{n} = \theta + \delta_{2} n^{-1/2}$ throughout the proof. Under \textbf{Assumption 1}, we have the following Taylor expansion:
\[
f_{n} = f_{0} + (\xi_{n}-\xi_{0})^{2} v_{\xi} + (\theta_{n}-\theta)^\top v_{\theta} + r_{n}(x_{i}; \xi_{n}, \theta_{n}). 
\]
Denoting $w_{i} = 2 ( f_{n}/f_{0}-1)$, we have 
\[
w_{i} = 2 (\xi_{n}-\xi_{0})^{2} \frac{v_{\xi}}{f_{0}} + 2 (\theta_{n}-\theta)^\top \frac{v_{\theta}}{f_{0}} + 2 R_{i}. 
\]
To show that under the modified DQM condition, the log-likelihood ratio admits a quadratic approximation, we use results in \textbf{Lemma 1}. \\

The log-likelihood ratio can be represented as 
\begin{align*}
\Lambda_{n} & =\sum_{i} \log \frac{p(x_{i} ; \xi_{n}, \theta_{n})}{p(x_{i} ; \xi_{0}, \theta)} = \sum_{i} 2 \log \frac{f_{n}}{f_{0}} = \sum_{i} 2 \log (1+ w_{i}/2)\\
&=\sum_{i} w_{i} - \frac{1}{4} \sum_{i}w_{i}^{2} + \frac{1}{2} \sum_{i}w_{i}^{2} \beta (w_{i})
\end{align*}
with $\beta(x) \to 0$ as $x\to 0$. \\

Using (\ref{eq:s2}) in \textbf{Lemma 1} and with $S_{n}=(S_{\xi, n}, S_{\theta,n}^\top)^\top$ and $J$ defined in \textbf{Theorem \ref{theorem1}}, we have 
\[
\sum_{i}w_{i} = 2 \frac{\delta_{1}^{2}}{\sqrt{n}} \sum_{i} \frac{v_{\xi}}{f_{0}} + 2 \frac{\delta_{2}^\top}{\sqrt{n}} \sum_{i} \frac{v_{\theta}}{f_{0}} + 2 \sum_{i} R_{i}
= t^\top S_{n} - \frac{1}{4} t ^\top J t + o_{P}(1)
\] 

Using (\ref{eq:s1}) and (\ref{eq:s3}) in \textbf{Lemma 1}, we have 
\[
\begin{array}{ll}
\sum_{i}w_{i}^{2}& = \sum_{i} \left (\frac{2\delta_{1}^{2}}{\sqrt{n}} \frac{v_{\xi}}{f_{0}} + \frac{2\delta_{2} ^\top}{\sqrt{n}} \frac{v_{\theta}}{f_{0}} + 2R_{i} \right)^{2}\\
&= t^\top J t + o_{P}(1) + 4 \sum_{i}R_{i}^{2} + 4 \sum_{i} R_{i} \left (\frac{2\delta_{1}^{2}}{\sqrt{n}} \frac{v_{\xi}}{f_{0}} + \frac{2\delta_{2}^\top}{\sqrt{n}} \frac{v_{\theta}}{f_{0}} \right) \\
&= t^\top J t + o_{P}(1)
\end{array}
\]\\

Lastly, we need to show that $\sum_{i} w_{i}^{2} \beta (w_{i}) = o_{p}(1)$. First note that using (\ref{eq:s4}) and (\ref{eq:s5}) in \textbf{Lemma 1}, we have \\
\[
\begin{array}{cc}
\mathbb{P} \left (\underset{1 \leq i \leq n}{\max} |w_{i}| > \epsilon \right) &\leq  \delta_{1}^{2} \mathbb{P} \left (\underset{1 \leq i \leq n}{\max} \left |\frac{2}{\sqrt{n}} \frac{v_{\xi}}{f_{0}} \right| >\epsilon \right) + \delta_{2}^\top \mathbb{P} \left (\underset{1 \leq i \leq n}{\max}  \left| \frac{2}{\sqrt{n}} \frac{v_{\theta}}{f_{0}} \right| > \epsilon \right) \\
&+ 2 \mathbb{P} \left (\underset{1 \leq i \leq n}{\max} |R_{i}| > \epsilon \right) \to 0
\end{array}
\]
Since when $w_{i} \to 0$, $\beta(w_{i})\to 0$, we have $\underset{1 \leq i \leq n}{\max}|\beta(w_{i})| = o_{p}(1)$. By H\"older's inequality, 
\[
\sum_{i} w_{i}^{2}\beta(w_{i}) \leq \underset{1 \leq i \leq n}{\max} |\beta (w_{i})| \sum_{i}w_{i}^{2} = o_{P}(1) O_{P}(1) = o_{P}(1). 
\]

Therefore, the log-likelihood ratio is approximated by 
\[
\begin{array}{ll}
\Lambda_{n} &= \sum_{i}w_{i} - \frac{1}{4} \sum_{i}w_{i}^{2} + \frac{1}{2} \sum_{i}w_{i}^{2} \beta(w_{i}) \\
& = t ^\top S_{n} - \frac{1}{4} t^\top J t -\frac{1}{4}t^\top J t + o_{p}(1)\\
&= t^\top S_{n} - \frac{1}{2} t^\top J t +o_{P}(1)
\end{array}
\]
\hfill{\rule{2mm}{2mm}}\\

\textbf{Proof of Corollary 1} 
Since $S_{n}$ is a normed iid sum, by the central limit theorem,
\[
S_{n} \overset{P_{n, \xi_{0}, \theta}}{\leadsto} \mathcal{N}(0, J)
\]
The zero asymptotic mean of $S_{n}$ is provided by (\ref{eq:mean0}) in \textbf{Lemma 1}, then the asymptotic variance for $S_{n}$ is $J$ as defined in Theorem 1. \\

The quadratic approximation for $\Lambda_{n}$ established in \textbf{Theorem 1} together with the joint normality of $S_{n}$ leads to the LAN property of the sequence of model $P_{n, \xi_{n}, \theta_{n}}$. Furthermore, we have 
\[
\Lambda_{n} \overset{P_{n, \xi_{0}, \theta}}{\leadsto} \mathcal{N}(-\frac{1}{2}t^\top J t, t^\top J t).
\]
By LeCam's first lemma (see e.g. \citeasnoun[Lemma 6.4]{vdv98}), $P_{n, \xi_{n}, \theta_{n}}$ and $P_{n, \xi_{0}, \theta}$ are mutually contiguous. \hfill{\rule{2mm}{2mm}}\\

\textbf{Proof of Theorem \ref{theorem2}}
The sequence of experiments $\mathcal{E}_{n}$ converges to a shifted Gaussian $\mathcal{N}(t, J^{-1})$ as a result of Theorem 9.4 in \citeasnoun{vdv98}. The log-likelihood ratio process of observing one sample from $\mathcal{N}(t, J^{-1})$ is 
\[
\log \frac{d\mathcal{N}(t, J^{-1})}{d\mathcal{N}(\textbf{0}, J^{-1})} (Y) = t^\top J Y - \frac{1}{2} t^\top J t
\]
It suffices to show that $J^{-1}S_{n}$ converges to the distribution of $Y$ under the null. \textbf{Corollary \ref{cor1}} establishes $S_{n} \overset{P_{n, \xi_{0}, \theta}}{\leadsto} \mathcal{N} \left (\textbf{0}, J \right )$, we thus have $J^{-1}S_{n} \overset{P_{n, \xi_{0}, \theta}}{\leadsto} \mathcal{N} \left (\textbf{0}, J^{-1} \right)$. \\

The optimal test statistic for $H_{0}: \delta_{1}=0$ against $H_{a}: \delta_{1} \neq 0$ in the limit experiment is the first element in $Y$. The sequence of test statistics from the original experiment $\mathcal{E}_{n}$ that matches with the first element in $Y$ is the $C(\alpha)$ statistic, 
\[
Z_{n} = (J_{\xi \xi} - J_{\xi \theta} J_{\theta \theta}^{-1} J_{\theta \xi})^{-1/2} (S_{\xi, n} - J_{\xi \theta}J_{\theta \theta}^{-1} S_{\theta, n}).
\]
Notice the rescaling in $Z_{n}$ is needed to obtain a unit asymptotic variance for the test statistic. \hfill{\rule{2mm}{2mm}}\\

\textbf{Proof of Corollary 2} 
Since $\xi$ is a scalar and $S_{n} \overset{P_{n, \xi_{0}, \theta}}{\leadsto} \mathcal{N}(0, J)$ under $H_{0}$, it is immediate that the asymptotic null distribution for $Z_{n}$ is $\mathcal{N}(0, 1)$. \\

We can now use LeCam's third lemma (see e.g. \citeasnoun[Example 6.7]{vdv98}) to derive the asymptotic distribution for $Z_{n}$ under local alternatives. We are interested in the local alternative that $\xi_{n}=\xi_{0}+\delta_{1}n^{-1/4}$ and nuisance parameter $\theta$ is left unspecified as in the null, hence we set $\delta_{2}=0$ in the log-likelihood ratio expansion. Under $H_{0}$, 
\[
(Z_{n}, \Lambda_{n}) \overset{P_{n, \xi_{0}, \theta}}{\leadsto } \mathcal{N} \begin{pmatrix}
\begin{pmatrix}
0\\ 
-\frac{1}{2} \delta_{1}^{4} J_{\xi \xi}
\end{pmatrix} & \begin{pmatrix}
1 & \sigma_{12}\\ 
\sigma_{12} & \delta_{1}^{4} J_{\xi \xi}
\end{pmatrix}
\end{pmatrix}
\]
with $\sigma_{12} = \Cov (Z_{n}, \Lambda_{n}) = \delta_{1}^{2} (J_{\xi \xi}-J_{\xi \theta}J_{\theta \theta}^{-1}J_{\theta \xi})^{1/2}$. With $\delta_{2}=0$, \textbf{Corollary 1} implies that $P_{n, \xi_{n}, \theta}$ are mutually contiguous to $P_{n, \xi_{0}, \theta}$, then LeCam's third lemma implies,
\[
Z_{n} \overset{P_{n, \xi_{n}, \theta}}{\leadsto} \mathcal{N}(\sigma_{12}, 1).  
\] 
\hfill{\rule{2mm}{2mm}} 

\textbf{Proof of Theorem 3}
Define the class of functions: 
\[
\FF_n := \Big\{x\mapsto (g(x,\theta) - g(x,\eta))\Big| \|\theta - \eta\| \leq \delta_n \Big\}.
\]
If $\hat \theta$ is a $\sqrt{n}$-consistent estimator of $\theta$, and $\delta_n = O(n^{-k})$ with $k <1/2$, we obtain that with probability tending to one
\[
\Big|Z_n(\hat\theta) - Z_n(\theta)\Big| \leq \sup_{f \in \FF_n} |\mathbb{G}_n(f)|
\]
where $\mathbb{G}_n(f) := n^{-1/2}\sum_i (f(X_i) - \E f(X_i))$ denotes the empirical process indexed by $\FF_n$. Proving $Z_n(\hat \theta) - Z_n(\theta) = o_P(1)$ thus amounts to establishing asymptotic equicontinuity of the process $\mathbb{G}_n$ with respect to the Euclidean norm. \\

Let the parameter space near true $\theta$, $U_{\delta_{n}}(\theta)$, be covered by balls with radius $\eps^{1/\gamma}$, the number of balls can be upper bounded by $C_{1}\eps^{-p/\gamma}$ with $C_{1}$ as a constant that does not depend on $n$ and $p$ being the dimension of the nuisance parameter space. Then for $\forall \eta \in U_{\delta_{n}}(\theta)$, $\exists N_{\eta}$, such that 
\[
\| \eta - \eta_{N_{\eta}} \| \leq \eps^{1/\gamma}
\]
The condition on $g$ in \textbf{Assumption \ref{assumption2}} implies
\[
|g(x,  \eta)-g(x,  \eta_{N_{\eta}})| \leq \|\eta - \eta_{N_{\eta}} \|^{\gamma}H(x) \leq \eps H(x)
\]
It follows that the bracketing number, $N_{[\ ]}(\eps \|H\|_{2}, \FF_n, \mathcal{L}_{2}(P_{n, \xi_n, \theta}))$ is bounded from above by $C_{2} \eps^{-p/\gamma}$. \\

Furthermore, the assumption also implies that for $f \in \FF_n$, $\|f\|_{P_{n}, 2} \leq \delta_{n}^{\gamma} \|H\|_{P_{n}, 2}$ with $\mathcal{L}_{2}(P_{n, \xi_{n}, \theta})$-norm. We can now apply Theorem 2.14.2 in \citeasnoun{vandwell1996} and get 
\[
\mathbb{E}_{P_{n, \xi_{n}, \theta}}(\sup_{f \in \FF_n} |\mathbb{G}_n(f)|) \leq J_{[\ ]}(\delta_{n}^{\gamma}, \FF_{n}, \mathcal{L}_{2}(P_{n, \xi_{n}, \theta})) \|H\|_{P_{n}, 2} + \sqrt{n} \mathbb{E}_{P_{n, \xi_{n}, \theta}}[H(X)I\{H(X)>\sqrt{n}a(\delta_{n}^{\gamma})\}]
\]
where the bracketing integral is defined as
\[
J_{[\ ]}(\delta_{n}^{\gamma}, \FF_{n}, \mathcal{L}_{2}(P_{n, \xi_{n}, \theta}))=
 \int _{0}^{\delta_{n}^{\gamma}} \sqrt{1+\log N_{[\ ]}(\eps \|H\|_{P_{n}, 2}, \FF_n, \mathcal{L}_{2}(P_{n, \xi_{n}, \theta}))} d\eps
 \]
  and 
 \[
  a(\delta_{n}^{\gamma}) = \delta_{n}^{\gamma} \|H\|_{2}/\sqrt{1+ \log N_{[\ ]}(\delta_{n}^{\gamma} \|H\|_{P_{n}, 2}, \FF_n, \mathcal{L}_{2}(P_{n, \xi_{n}, \theta}))}.\\
\]
  
Provided that $\delta_{n} \to 0$, we have for $n$ large enough, 
\[
J_{[\ ]}(\delta_{n}^{\gamma}, \FF_{n}, \mathcal{L}_{2}(P_{n, \xi_{n}, \theta})) \leq \int _{0}^{\delta_{n}^{\gamma}} 1+ \log (C_{2} \eps^{-p/\gamma}) d\eps \to 0
\]
Since $H(x)$ is square integrable for all $n$ by \textbf{Assumption \ref{assumption2}}, the first term goes to zero. \\

The upper bound for the bracketing number also yields a lower bound for $a(\delta_{n}^{\gamma})$ that is for $\delta_{n}$ sufficiently small, 
\[
a(\delta_{n}^{\gamma}) \geq \frac{\delta_{n}^{\gamma} \|H\|_{P_{n}, 2}}{\sqrt{1+\log (C_{2} \delta_{n}^{-p})}} := k_{n} \to 0
\]
As long as $k_{n}$ converges to zero slower than $c_{n}$, \textbf{Assumption \ref{assumption2}} ensures that the second term also tends to zero. \\

The last step is to check that $\underset{f \in \FF_n}{\sup} \frac{1}{\sqrt{n}} \sum_{i} \mathbb{E}_{P_{n, \xi_{n}, \theta}}(f(X_{i})) = o(1)$ so that $\underset{f\in \FF_n}{\sup}|\mathbb{G}_{n}(f)|$ is the correct upper bound. This is trivially true under the null, where $\xi_{n}=\xi_{0}$ for all $n \in \mathbb{N}$, since $\mathbb{E}_{P_{n, \xi_{0}, \theta}}(g(X_{i}, \theta)) =  \mathbb{E}_{P_{n, \xi_{0}, \theta}}(g(X_{i}, \hat \theta))=0$. Under local alternatives with $\xi_{n} = \xi_{0}+ \delta_{1} n^{-1/4}$ and given the i.i.d. assumption on the sample, it suffices to show that 
\[
\underset{\|\eta - \theta\|\leq \delta_n}{\sup} \sqrt{n} \int (g(x, \eta)-g(x, \theta)) p(x; \xi_n, \theta)dx = o(1)
\]
Denote $p_n=p(x; \xi_n, \theta)$ and $p_0=p(x; \xi_0, \theta)$, we have the following expansion 
\[
\begin{array}{lrl}
& \sqrt{n} \int (g(x, \eta)-g(x, \theta)) p_{n} dx  \\
& = &  \sqrt{n} \int \Big ((g(x, \eta)-g(x, \theta)) (\sqrt{p_0}+(\xi_n-\xi_0)^{2} v_{\xi}(x) + r_{n} \Big)\sqrt{p_n}dx \\
& =&   \sqrt{n} \int (g(x, \eta)- g(x, \theta))   \sqrt{p_0}\sqrt{p_n} dx\\
& & {} +  \sqrt{n}(\xi_n-\xi_0)^2 \int (g(x, \eta)-g(x, \theta)) \sqrt{p_n}v_{\xi}(x)dx \\
& &  {}+  \sqrt{n} \int (g(x, \eta)-g(x, \theta)) \sqrt{p_n}r_n dx
\end{array}
\]
The last two terms are $o(1)$ uniformly over $\eta$ for $\|\eta-\theta\| \leq \delta_n$ due to the DQM condition in \textbf{Definition 1} and assumption on $g$ in \textbf{Assumption \ref{assumption2}}. Since Cauchy-Schwarz inequality implies that with respect to $\mathcal{L}_{2}(\mu)$-norm, 
\[
\begin{array}{ll}
|\int (g(x, \eta)-g(x, \theta)) \sqrt{p_n} v_{\xi}(x)dx| &\leq \| (g(x, \eta)-g(x, \theta))\sqrt{p_0}\|_{\mu, 2} \|v_{\xi}\|  _{\mu, 2}\\
& \leq \|\eta-\theta\|^{\gamma} \|H\|_{P_{n}, 2}\|v_{\xi}\|_{\mu, 2} =o(1).
\end{array}
\]
Similarly, \\
\[
\begin{array}{ll}
|\sqrt{n} \int (g(x, \eta)-g(x, \theta)) \sqrt{p_n} r_{n} dx| &\leq \| (g(x, \eta)-g(x, \theta))\sqrt{p_0}\|_{\mu, 2} \sqrt{n}\|r_n\|_{\mu, 2}  
=o(1).
\end{array}
\]
The first term is also $o(1)$ by expanding $\sqrt{p_n}$ again and applying Cauchy-Schwarz inequality in a similar fashion.   \hfill{\rule{2mm}{2mm}} \\

\textbf{Proof of Theorem \ref{theorem4}}
As in the proof of \textbf{Theorem \ref{theorem2}}, the limit of the sequence $\upsilon_{n}$ is a shifted Gaussian experiment $Y \sim \mathcal{N}(t, J^{-1})$ but now with $t^\top=(\delta_{1}^{2}, \delta_{2}^{2},2\delta_{1}\delta_{2}, \delta_{3}^\top)$. An equivalent limit experiment observes $X \sim \mathcal{N}( Jt, J)$ with $X=JY$, because the likelihood ratio process of $\frac{d\mathcal{N}(t, J^{-1})}{d \mathcal{N}(0, J^{-1})}(Y)$ is identical to that of $\frac{d\mathcal{N}(t^\top J, J)}{d \mathcal{N}(0, J)}(X)$.  \\

To be more explicit, denoting the first three elements of $X$ to be $X_{\xi}$, and the rest to be $X_{\theta}$, we have under the alternative, 
\[
\begin{pmatrix} 
X_{\xi} \\
X_{\theta}
\end{pmatrix}
\overset{\mathcal{D}}{=} \mathcal{N} \left ( \begin{pmatrix} 
J_{\xi \xi} & J_{\xi \theta}\\
J_{\theta \xi} & J_{\theta \theta}
\end{pmatrix} \begin{pmatrix} 
t_{\xi} \\
t_{\theta}
\end{pmatrix}, J\right)
\]
with $t_{\xi}=(\delta_{1}^{2}, \delta_{2}^{2}, 2\delta_{1}\delta_{2})^\top$ and $t_{\theta}=\delta_{3}^\top$. \\

To focus on testing for zero restrictions on $t_{\xi}$, we find the conditional distribution of $X_{\xi}$ on $X_{\theta}$ to be 
\[
\tilde X_{\xi} = X_{\xi}-J_{\xi \theta}J_{\theta \theta}^{-1} X_{\theta} \overset{\mathcal{D}}{=} \mathcal{N}((J_{\xi \xi}-J_{\xi \theta}J_{\theta \theta}^{-1} J_{\theta \xi})t_{\xi}, J_{\xi \xi}-J_{\xi \theta}J_{\theta \theta}^{-1} J_{\theta \xi}).
\]
The matched statistic from the original experiment is then 
\[
\tilde S_{\xi, n} = S_{\xi, n} - J_{\xi \theta}J_{\theta \theta}^{-1} S_{\theta, n}
\]
Under $H_{0}$,  $\tilde S_{\xi, n} $ follows $\mathcal{N}(0, \Sigma)$ with $\Sigma = J_{\xi \xi}-J_{\xi \theta}J_{\theta \theta}^{-1} J_{\theta \xi}$, and under local alternative, its asymptotic distribution is $\mathcal{N}(\Sigma t_{\xi}, \Sigma)$. \\

Notice we can decompose $\tilde S_{\xi,n} \Sigma^{-1} \tilde S_{\xi,n}$ into two independent pieces as $u_{n}^\top \Sigma_{11.2} u_{n} + w_{3n}^\top w_{3n}$. Let the Cholesky decomposition of $\Sigma_{11.2}$ be such that $\Lambda \Lambda^\top = \Sigma_{11.2}$, then
$w_{n} := \Lambda^{-1} u_{n}  \overset{\mathbb{P}_{n, \xi_{0}, \theta}}{\leadsto} \mathcal{N}(0, I)$ and 
$w_{n} \overset{\mathbb{P}_{n, \xi_{n}, \theta_{n}}}{\leadsto} \mathcal{N}(\Lambda^\top \begin{pmatrix} \delta_{1}^{2} \\ \delta_{2}^{2} \end{pmatrix}, I)$. Since $(\delta_{1}^{2}, \delta_{2}^{2}) \in \mathbb{R}_{+}^{2}$ and 
\[
\begin{pmatrix}
\eta_{1}\\
\eta_{2}
\end{pmatrix}
:=  \Lambda^\top \begin{pmatrix} \delta_{1}^{2}\\ \delta_{2}^{2} \end{pmatrix} = \begin{pmatrix}
\delta_{1}^{2} \sqrt{v_{1}} + \rho \delta_{2}^{2} \sqrt {v_{2}}\\
\sqrt{v_{2}} \sqrt{1-\rho^{2}} \delta_{2}^{2}
\end{pmatrix}
\]
The feasible parameter set is therefore the convex cone defined as, 
\[
\left \{ (\eta_{1}, \eta_{2}) \mid \eta_{2} \geq 0, \eta_{1}-\frac{\rho}{\sqrt{1-\rho^{2}}} \eta_{2} \geq 0  \right \}.
\]
For test statistic taking a value that falls outside of the feasible set, it needs to be projected onto the set. This yields the following four cases as illustrated in the figure. 

\begin{figure}
\includegraphics{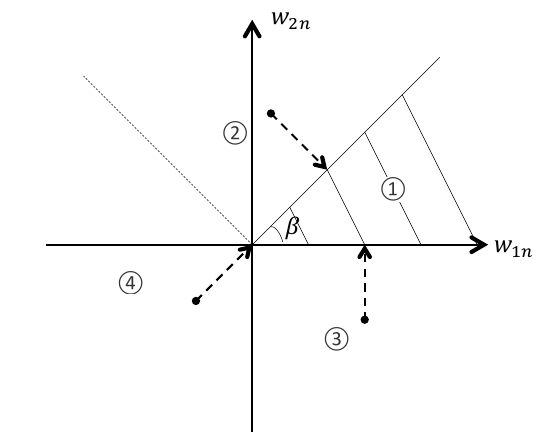}
\end{figure}

\textbf{Case 1}:
When the value of the test statistic $w_{n}$ falls into shaded area \circled{1}, the test statistics is the sum of squares of the elements of $w_n$ and $w_{3n}$ which are mutually independent:
\[
T_{n} = w_{1n}^{2} + w_{2n}^{2} + w_{3n}^{2} \sim \chi_{3}^{2}
\]

\textbf{Case 2}:
When the test statistic falls into area \circled{2}, we need to project $w_{n}$ onto the convex cone \circled{1}, which gives a point with coordinates $(\rho^{2}w_{1n} + \rho \sqrt{1-\rho^{2}}w_{2n}, \rho \sqrt{1-\rho^{2}}w_{1n}+(1-\rho^{2})w_{2n})$. The $C(\alpha)$ test statistic is hence:
\[
\begin{array}{ll}
T_{n} &= (\rho^{2}w_{1n} + \rho \sqrt{1-\rho^{2}}w_{2n})^{2}+ (\rho \sqrt{1-\rho^{2}}w_{1n}+(1-\rho^{2})w_{2n})^{2} + w_{3n}^{2}\\
&= (\rho w_{1n}+\sqrt{1-\rho^{2}}w_{2n})^{2} + w_{3n}^{2} \sim \chi_{2}^{2}
\end{array}
\]

\textbf{Case 3}:
When the test statistic $w_{n}$ falls in area \circled{3}, projecting onto the region \circled{1} yields $(w_{1n}, 0)$ and thus, 
\[
T_{n}=w_{1n}^{2} + w_{3n}^{2} \sim \chi_{2}^{2}
\]

\textbf{Case 4}:
Lastly, when $w_{n}$ falls into region \circled{4}, projecting onto region \circled{1} yields $(0,0)$ and hence,
\[
T_{n}=0+w_{3n}^{2} \sim \chi_{1}^{2}\\
\]

The asymptotic distribution of the $C(\alpha)$ test statistics is a mixture of $\chi^{2}$'s, for which the weights are characterized by the probability of falling into different regions. The angle $\beta$ spanned by the shaded area \circled{1} as marked in the figure is $\beta = \cos^{-1}(\rho)$, hence the probability of falling into region \circled{1} is $\frac{\beta}{2\pi}$. The probability of falling into \circled{2} and \circled{3} is $\frac{1}{2}$, leaves the probability of falling into \circled{4} as $(\frac{1}{2}-\frac{\beta}{2\pi})$. \hfill{\rule{2mm}{2mm}}

\section{Computational details in examples}

\subsection{Joint test for Gaussian panel data model}
The information matrix for $(\xi, \theta)=(\xi_{1}, \xi_{2}, \mu_{0}, \sigma_{0}^{2})$ is 
\[
I = \begin{pmatrix}
I_{\xi \xi} & I_{\xi \theta}\\
I_{\theta \xi} & I_{\theta \theta}
\end{pmatrix} = 
\frac{NT}{\sigma_{0}^{4}} \begin{pmatrix}
2T & \sigma_{0}^{2} & 0 & 1 \\
\sigma_{0}^{2} & (T+3)\sigma_{0}^{4}/2 & 0 & \sigma_{0}^{2}/2 \\
0 & 0 & \sigma_{0}^{2} & 0\\
1 &\sigma_{0}^{2}/2 & 0 & 1/2
\end{pmatrix}
\]
We further find 
\[
I_{\xi.\theta} = I_{\xi \xi}-I_{\xi \theta}I_{\theta \theta}^{-1}I_{\theta \xi} =
\begin{pmatrix}
2NT(T-1) / \sigma_{0}^{4}& 0 \\
0 & NT(T/2+1)
\end{pmatrix}
\]
and 
\[
I_{\xi \theta} I_{\theta \theta}^{-1} = \begin{pmatrix} 0 & 2 \\ 0 & \sigma_{0}^{2} \end{pmatrix}. 
\]\\

As we have remarked in Section \ref{2.3.3}, the diagonality of $I_{\xi.\theta}$ provides much convenience for finding the optimal test statistics. Denote 
\[
T_{n} := \begin{pmatrix} t_{1n}\\t_{2n} \end{pmatrix}
= I_{\xi.\theta}^{-1/2} \begin{pmatrix}
\sum_{i}v_{i1} - 2 \sum_{i} v_{4i} \\
\sum_{i}v_{2i} - \sigma_{0}^{2} \sum_{i} v_{4i}
\end{pmatrix} 
 = \begin{pmatrix}
(2NT(T-1)/\sigma_{0}^{4})^{-1/2} \left (\sum_{i} (\frac{\bar y_{i.}-\mu_{0}}{\sigma_{0}^{2}/T})^{2} - NT/\sigma_{0}^{2} \right) \\
(NT(T/2+1))^{-1/2} \left ( \sum_{i} (Z_{i} - T/2)^{2}-NT/2 \right )
\end{pmatrix}
\]
Replacing $(\mu_{0}, \sigma_{0}^{2})$ by their MLEs yields the joint $C(\alpha)$ test. 

\section{Claim in Section 4}
Here we provide the detail derivation for the claim in Section 4 that the reparameterization adopted in \citeasnoun{Chesher84} and \citeasnoun{Cox83} for heterogeneity test requires extra moment conditions on $U$ for second derivative of log density with respect to the test parameter to be bounded. 

\begin{proposition}
For iid random variable $Y_{1}, \dots, Y_{n}$ each with density function $\int p(y; \lambda_{0} + \tau \sqrt{\eta} u_{i})dF(u_{i})$, where $U_{i}$ is a random variable with zero mean and unit variance. The second-order derivative of the log density with respect to $\eta$ evaluated under $\eta=0$ is unbounded unless $\mathbb{E}(U^{3})=0$ and $\mathbb{E}(U^4)<\infty$. 
\end{proposition}
\textbf{Proof}
Denote the log density as $l=\log \int p(y; \lambda_{0}+\tau \sqrt{\eta} u_{i})dF(u_{i})$. The first order derivative with respect to $\eta$ is
\[
\left. \nabla_{\eta}l \right |_{\eta=0} = \frac{\tau \int \nabla_{\lambda}p(y; \lambda_{0}) u dF(u)}{2 \sqrt{\eta} \int p(y; \lambda_{0})dF(u)}
=\frac{\tau^{2}}{2} \mathbb{E}(U^{2}) \frac{\nabla_{\lambda}^{2} p(y; \lambda_{0})}{p(y; \lambda_{0})}
\]
The last step is obtained by applying the l'H\^opital's rule. \\

The second order derivative is 
\[
\begin{array}{ll}
\left. \nabla_{\eta} ^{2} l \right |_{\eta=0} &= \left. \frac{\tau^{2} \sqrt{\eta} \int \nabla_{\lambda}^{2}p(y; \lambda_{0})u^{2}dF(u) - \tau \int \nabla_{\lambda}p(y; \lambda_{0})udF(u)}{4\eta \sqrt{\eta} \int p(y; \lambda_{0})dF(u)} \right |_{\eta=0}- \left (\left. \nabla_{\eta}  l \right |_{\eta=0} \right)^{2}\\
&= \left. \frac{\tau^{3} \int \nabla_{\lambda}^{3} p(y; \lambda_{0})u^{3}dF(u)}{12 \sqrt{\eta} \int p(y; \lambda_{0})dF(u)} \right|_{\eta=0} - \left (\left. \nabla_{\eta}  l \right |_{\eta=0} \right)^{2}
\end{array}
\]
Provided that $\nabla_{\lambda}^{3} p(y; \lambda_{0})$ is not degenerately zero,  $\nabla_{\eta}^{2}l$ is unbounded unless $\mathbb{E}(U^3)=0$ and $\mathbb{E}(U^{4})<\infty$ so that we can apply l'H\^opital's rule again and get 
\[
\left. \nabla_{\eta}^{2}l\right |_{\eta=0} = \frac{\tau^{4}}{12} \left [ \mathbb{E}(U^4) \frac{\nabla_{\lambda}^{4} p(y; \lambda_{0})}{p(y; \lambda_{0})} - 3 \mathbb{E}(U^2)^2\Big (\frac{\nabla_{\lambda}^{2}p(y; \lambda_{0})}{p(y; \lambda_{0})}\Big )^2 \right] < \infty
\] \hfill{\rule{2mm}{2mm}}

\end{appendix}
\end{document}